\documentclass{article}

\usepackage{amsthm,amsfonts,amssymb}
\usepackage{amsmath}
\usepackage[all]{xy}

\newtheorem{thm}{Th\'eor\`eme}[section]
\newtheorem{cor}[thm]{Corollaire} 
\newtheorem{lem}[thm]{Lemme}
\newtheorem{prop}[thm]{Proposition} 

\theoremstyle{definition} 
\newtheorem{defn}[thm]{D\'efinition}
\theoremstyle{definition} 
\newtheorem{ex}[thm]{Exemple}
\theoremstyle{definition}

\newcommand{\R}{\mathbb R}
\newcommand{\C}{\mathbb C} 
\newcommand{\Z}{\mathbb Z}
\newcommand{\N}{\mathbb N}

\newcommand{\F}{\mathcal F}
\newcommand{\K}{\mathcal K} 
\renewcommand{\L}{\mathcal L}
\renewcommand{\H}{\mathcal H} 

\newcommand{\mois}{%
  \ifcase\month\or Janvier\or F\'evrier\or Mars\or Avril\or Mai\or
  Juin\or Juillet\or Ao\^ut\or Septembre\or Octobre\or Novembre\or
  D\'ecembre\fi }

\title{{\bf La caract\'eristique d'Euler des feuilletages mesur\'es}}
\date{}

\author{Miguel Berm\'udez\\ 
\\
{\small Institut Camille Jordan, Universit\'e Claude
Bernard - Lyon 1} \\ 
{\small 43, Blvd. du 11 novembre 1918, F-69622 Villeurbanne
(France)} \\ 
{\small e-mail: bermudez@igd.univ-lyon1.fr}}

\begin{document}
\maketitle

\renewcommand{\abstractname}{R\'esum\'e}

\section*{Introduction}
Soit $\mathcal F$ un feuilletage orientable de dimension $n$ et classe
$C^2$ sur une d'une vari\'et\'e diff\'erentiable compacte $M$, munie
d'une mesure transverse invariante $\mu$. On appellera dans ce papier
le triplet $(M,\mathcal F,\mu)$ un {\em feuilletage mesur\'e}. En
int\'egrant les $p$-formes diff\'erentielles de $M$ sur les
feuilles puis relativement \`a la mesure transverse on d\'efinit un
courant ferm\'e $C_\mu:\Omega^n(M)\to \R$ appel\'e {\em courant de
  Ruelle-Sullivan} \cite{RS}. Sa classe d'homologie $[C_\mu]\in
H_n(M;\R)$ est souvent interpr\'et\'ee comme la classe fondamentale du
feuilletage mesur\'e $(M,\mathcal F,\mu)$. On d\'efinit la {\em
  caract\'eristique d'Euler} de $(M,\mathcal F,\mu)$ par $$
\chi(M,\mathcal F,\mu)=\langle e(T\mathcal F),[C_\mu]\rangle $$
o\`u $e(T\mathcal F)$ est la classe d'Euler du $n$-fibr\'e tangent au
feuilletage  $T\mathcal F\to M$. Dans ce papier on
fait le lien entre la caract\'eristique d'Euler et les sections du
fibr\'e tangent $T\mathcal F$, et retrouvons dans ce cadre des
versions des r\'esultats classiques bien connus pour les vari\'et\'es
compactes. Le contexte o\`u nos r\'esultats sont valables est en fait
bien plus large que celui des feuilletages des vari\'et\'es compactes,
et contient celui 
des laminations (l'espace ambiant n'est plus un vari\'et\'e) et des
feuilletages mesurables (la r\'egularit\'e transverse est suppos\'ee
seulement mesurable). Par souci de clart\'e nous avons voulu nous
restreindre dans ce papier au cas plus simple des feuilletages des
vari\'et\'es compactes. Tous les objets seront donc topologiques, \`a 
l'exception des champs 
de vecteurs et les m\'etriques de Riemann, qui seront suppos\'es
continus le long des feuilles mais seulement mesurables
transversalement. De tels objets seront souvent appel\'es {\em de
  classe $MC^0$} pour mesurable et continu de classe
$C^0$. L'int\'er\^et de consid\'erer des objets dans cette cat\'egorie
est illustr\'e par les \'equivalences \'enonc\'ees dans les
th\'eor\`emes $B$ et $C$, en particulier le second. Ils ne sont pas
vrais si l'on consid\`ere des champs transversalement continus !  

Nous d\'emontrons dans ce travail une s\'erie de th\'eor\`emes qui
clarifient la signification g\'eom\'etrique et ergodique de la
caract\'eristique d'Euler, en reliant celle-ci \`a l'existence de
champs tangents sans z\'ero de classe $MC^0$ sur le feuilletage. On
commence par une version feuillet\'ee du th\'eor\`eme bien connu de
Poincar\'e-Hopf: 

\newtheorem*{thmA}{Th\'eor\`eme A}
\begin{thmA}\label{thm:A} {\em Soit $\bf x$ un champ tangent
    \`a $\mathcal F$ de classe $MC^0$ et soit $O_{\bf x}$ l'ensemble
    des z\'eros de $\bf x$. Si la trace 
    de $O_{\bf x}$ sur $\mu$-presque toute feuille est discr\`ete
    dans cette feuille, i.e. si $O_{\bf x}$ est une transversale
    mesurable de $(M,\mathcal F,\mu)$, alors l'indice local $ind_{\bf
      x}(x)$ est d\'efini pour $\mu$-presque tout $x\in O_{\bf x}$.
    Si la fonction mesurable \`a valeurs enti\`eres $ind_{\bf x}$ est
    dans $L^1(O_{\bf x},\mu)$ alors on a: 
    $$
    \chi(M,\mathcal F,\mu)= \int_{O_{\bf x}} ind_{\bf x}(x)d\mu(x).
    $$
  }
\end{thmA}
Il avait d\'ej\`a \'et\'e remarqu\'e par Connes dans \cite{Con2} que ce
r\'esultat \'etait vrai pour un champ transversalement continu et
transverse \`a la section nulle de $T\mathcal F$. Ce fait
est facile \`a prouver et d\'ecoule automatiquement de la construction
topologique de la classe $e(T\mathcal F)$. L'id\'ee est de voir la
fonction indice $ind_{\bf x}$ comme un $n$-cocycle cellulaire de $M$,
puis de montrer que sa classe de cohomologie ne d\'epend pas du champ
choisi. On d\'efinit ainsi un \'el\'ement de $H^n(M)$ qui s'av\`ere
\^etre $e(T\mathcal F)$ pour des raisons de naturalit\'e. Ce
raisonnement ne s'applique \'evidemment pas au cas des champs
transversalement mesurables. Ils ne d\'efinissent pas une classe dans
$H^n(M)$, mais un \'el\'ement dans un groupe de cohomologie simplicial
mesurable similaire \`a ceux utilis\'es dans \cite{MS}. Celui-ci
ne d\'epend pas de la topologie de $M$ mais seulement de celle des
feuilles, et puis aussi de la dynamique mesurable transverse du
feuilletage. Le th\'eor\`eme A ne se ram\`ene donc pas \`a \cite{Con2}.

On d\'eveloppera ici les outils de 
base n\'ecessaires pour expliciter l'information
g\'eom\'etrique et ergodique contenue dans ces groupes de cohomologie
mesurable, et pr\'eciser la relation qui existe entre cette
cohomologie et la cohomologie classique de $M$ o\`u l'on situe
habituellement les invariants qui nous int\'eressent.

\newtheorem*{thmB}{Th\'eor\`eme B}
\begin{thmB}\label{thm:B} 
{\em 
  Soit $(M,\mathcal F,\mu)$ un feuilletage mesur\'e ergodique.  Alors
  les deux conditions suivantes sont \'equivalentes: 
  \begin{enumerate}
  \item $\chi(M,\mathcal F,\mu)=0$.
  \item Pour tout $\epsilon >0$ il existe un champ tangent $\mathbf x$
    de classe $MC^0$ et \`a singularit\'es non d\'eg\'en\'er\'ees tel
    que $\mu(O_{\bf x})<\epsilon$.   
  \end{enumerate}
}
\end{thmB}

Rappelons qu'un feuilletage mesur\'e $(M,\mathcal
F,\mu)$ est dit {\em moyennable} s'il existe une famille de
fonctionnelles continues (dites moyennes) $m_x:L^\infty(L_x)\to \R$ ($x\in
M$) v\'erifiant les conditions suivantes:
\begin{enumerate}
\item $m_x(f)\geq 0$ si $f\geq 0$;
\item $m_x(\mathbf 1)=1$;
\item $m_x=m_y$ si $x$ et $y$ appartiennent \`a la m\^eme feuille;
\item Pour toute fonction mesurable born\'ee $f:M\to \R$ la fonction
  $m(f):M\to \R$ d\'efinie par $m(f)(x)=m_x(f|L_x)$ est mesurable.
\end{enumerate}

Pour les feuilletages moyennables le th\'eor\`eme B peut \^etre
am\'elior\'e de la fa\c con suivante:

\newtheorem*{thmC}{Th\'eor\`eme C}
\begin{thmC}\label{thm:C}
  {\em Soit $(M,\mathcal F,\mu)$ un feuilletage mesur\'e ergodique
    moyennable.  Alors les deux conditions suivantes sont
    \'equivalentes: 
    \begin{enumerate}
    \item $\chi(M,\mathcal F,\mu)=0$;
    \item Le feuilletage poss\`ede un champ tangent sans z\'ero
      de classe $MC^0$.
    \end{enumerate}
}
\end{thmC}

Il existe des feuilletages non moyennables et \`a caract\'eristique
d'Euler nulle. Par exemple la suspension d'une repr\'esentation
fid\`ele du groupe fondamental d'une surface de genre $\geq 2$ dans
le groupe des isom\'etries d'une vari\'et\'e de Riemann compacte
fournit un feuilletage mesur\'e par plans hyperboliques non moyennable
(voir \cite{Zim}). Si on le multiplie par $S^1$, on obtient un
feuilletage mesur\'e de dimension trois, qui est
toujours non moyennable, mais  dont la caract\'eristique d'Euler est
nulle pour des raisons de dualit\'e. Avec un peu plus de travail on
peut aussi construire des feuilletages non moyennables et \`a
caract\'eristique d'Euler nulle en toute dimension $\geq 3$. En
dimension deux, par contre, la caract\'eristique d'Euler constitue un
invariant ergodique contenant une grande quantit\'e d'information
aussi bien g\'eom\'etrique que dynamique, comme le montre notre
prochain th\'eor\`eme.  

\'Etant donn\'ee une m\'etrique de Riemann de classe $MC^0$ sur un
feuilletage mesur\'e $(M,\mathcal F,\mu)$, on peut d\'efinir une mesure
globale sur l'espace $M$ qui est localement d\'efinie comme la mesure
produit du volume $vol_g$ sur les feuilles et de la mesure transverse
invariante $\mu$, et qu'on notera $\mu\otimes vol_g$. On dira
que $g$ est {\em \`a volume $\mu$-fini} si cette mesure est de masse
totale finie, i.e. $\mu\otimes vol_g(M)<\infty$. Remarquons enfin que
la compacit\'e de $M$ implique que toute m\'etrique de Riemann
continue sur le feuilletage mesur\'e $(M,\mathcal F,\mu)$ est \`a
g\'eom\'etrie born\'ee et volume $\mu$-fini. Ce n'est \'evidemment pas
le cas en g\'en\'eral des m\'etriques de classe $MC^0$.

\newtheorem*{thmE}{Th\'eor\`eme D}
\begin{thmE}\label{thm:D}
{\em 
Soit $(M,\mathcal F,\mu)$ un feuilletage mesur\'e ergodique de
dimension deux. Alors les sept conditions suivantes sont
\'equivalentes: 
\begin{enumerate}
\item $\chi(M,\mathcal F,\mu)=0$;
\item $(M,\mathcal F,\mu)$ poss\`ede un champ tangent sans z\'ero de
  classe $MC^0$;
\item Toute m\'etrique de Riemann de classe $MC^0$ \`a volume
  $\mu$-fini  et \`a g\'eom\'etrie born\'ee sur chaque feuille est
  parabolique sur $\mu$-presque toute feuille, i.e. le rev\^etement
  universel de $\mu$-presque toute 
  feuille est conform\'ement \'equivalent au plan euclidien $\C$;
\item $(M,\mathcal F,\mu)$ est moyennable et $\mu$-presque toutes ses
  feuilles sont des tores, des cylindres ou des plans;
\item Il existe une application mesurable $\rho:M\to \mathbb T^2$ qui
  est un rev\^etement en restriction \`a chaque feuille. Autrement
  dit, le feuilletage est isomorphe \`a la suspension d'une action
  mesurable ergodique de $\Z^2$ sur un espace de Lebesgue.
\item Le feuilletage est d\'efini par une action de $\R^2$
  sur $M$ de classe $MC^0$, i.e. il existe une action mesurable de
  $\R^2$ sur $M$ qui est continue et localement libre le long de
  chaque feuille. 
\item La feuilletage poss\`ede une m\'etrique de Riemann
  de classe $MC^0$ qui est plate et compl\`ete sur chaque
  feuille. 
\end{enumerate}
}
\end{thmE}

Les implications
(3)$\Rightarrow$(4)$\Rightarrow$(5)$\Rightarrow$(6)$\Rightarrow$(7)$\Rightarrow$(4)
sont prouv\'ees dans \cite{BH}. Pour d\'emontrer le th\'eor\`eme, on
\'etablira la cha\^ine d'implications
(6)$\Rightarrow$(2)$\Rightarrow$(1)$\Rightarrow$(3). La premi\`ere est
triviale et la deuxi\`eme est une cons\'equence du th\'eor\`eme
A. Nous consacrerons le paragraphe \S\ref{sec:preuve-du-theoreme-C}
\`a la preuve de l'implication (1)$\Rightarrow$(3). 

\medskip 
Ce papier comporte deux parties. La premi\`ere, compos\'ee
des \S1,\S2 et \S3, contient des d\'eveloppements g\'en\'eraux
n\'ecessaires pour la preuve des quatre th\'eor\`emes
\'enonc\'es. Tous les r\'esultats dans cette partie seront \'etablis
dans un cadre bor\'elien. On s'int\'eresse dans cette partie \`a des
applications {\em de classe $BC^0$}, i.e. continues le long des
feuilles et globalement bor\'eliennes. Au \S1 on d\'emontre quelques
r\'esultats basiques de la th\'eorie des 
feuilletages mesurables, dont le th\'eor\`eme d'approximation
simpliciale. Au \S2 on d\'eveloppe une th\'eorie de l'obstruction
analogue \`a celle de \cite{Eil} qui nous permet de r\'esoudre
compl\`etement le probl\`eme de l'extension d'une application de
classe $BC^0$. Au \S3 on applique les r\'esultats du \S2 au cadre des
fibr\'es localement triviaux pour donner une construction
g\'eom\'etrique des classes caract\'eristiques feuillet\'ees
mesurables. On y r\'esout le probl\`eme de la construction de sections
de classe $BC^0$ d'un fibr\'e localement trivial. La deuxi\`eme,
compos\'ee du \S4, contient la preuve des th\'eor\`emes A, B, C et D
proprement dite. On y pr\'ecise la notion de mesure transverse et
d'application {\em de classe $MC^0$}, i.e. une application continue le
long des feuilles qui est bor\'elienne une fois qu'on a enlev\'e un
ensemble de feuilles de mesure nulle.

\section{Quelques r\'esultats pr\'eliminaires}\label{sec:app-pri}
On d\'emontre ici quelques r\'esultats techniques basiques
n\'ecessaires pour la suite. Nous rappelons qu'un espace de Borel
standard est un espace 
mesurable isomorphe \`a un bor\'elien d'un espace de polonais. Presque
tous les espaces mesurables qui apparaissent en topologie sont de ce
type. Il est connu qu'un espace de Borel standard est soit discret
d\'enombrable, soit isomorphe \`a l'intervalle $[0,1]$ de la droite
r\'eelle. 

\subsection{Le th\'eor\`eme de Kallman}
\label{sec:Kallman}
La difficult\'e fondamentale qu'on retrouve au moment d'\'etablir les
r\'esultats \'enonc\'es dans l'introduction est la construction de
sections bor\'eliennes de certaines applications bor\'eliennes ou
continues. Le r\'esultat suivant, prouv\'e par Kallman, est un outil
fondamental pour r\'esoudre ces difficult\'es:

\begin{thm}[\cite{Kal}]\label{thm:thm_kallman}
  Soit $Y$ un espace polonais et $X$ un espace de Borel standard. Soit
  $f:Y\to X$ une application surjective bor\'elienne dont la fibre
  $f^{-1}(x)$ en chaque  point $x\in X$ est r\'eunion d\'enombrable de
  compacts de $Y$. Alors $f$ poss\`ede une section bor\'elienne,
  i.e. il existe une application bor\'elienne $s:X\to Y$ telle que
  $f\circ s(x)=x$ pour tout $x\in X$. 
\end{thm}

En fait on n'aura besoin, la plus part du temps, que du corollaire
suivant, dont une preuve peut \^etre trouv\'e aussi dans \cite{Kur}:

\begin{cor}
  Soit $f:Y\to X$ une application bor\'elienne surjective entre deux
  espaces de Borel standard. Si $f^{-1}(x)$ est d\'enombrable pour
  tout $x\in X$, alors $f$ poss\`ede une section bor\'elienne.
\end{cor}

\subsection{Familles bor\'eliennes d'applications}
\label{sec:familles_applications}
Soient $B$ et $F$ deux espaces polonais connexes. On consid\`ere
l'espace des applications continues $C(B,F)$ muni de la topologie
compacte-ouverte. Rappelons que cette topologie est d\'efinie par la
sous-base 
$$ 
\mathcal V(K,U)=\{f\in
C(B,F)~|~f(K)\subset V\} 
$$
o\`u $K$ parcourt les compacts de $B$ et
$U$ les ouverts de $F$. Si $B$ est compact alors l'espace $C(B,F)$ est
lui aussi un espace Polonais. En effet, si $d$ est une m\'etrique
compl\`ete sur $F$, alors la m\'etrique de la convergence uniforme
d\'efinie par $$
d^*(f,g)=\sup\{d(f(x),g(x))~|~x\in B\} $$
est
compl\`ete sur $C(B,F)$. En plus si $U_i$ est une base d\'enombrable
pour $F$, alors $U_i^*=\{f\in C(B,F)~|~f(B)\subset U_i\}$ est une base
d\'enombrable pour $C(B,F)$. Soit $T$ un espace de Borel standard et
$g:B\times T\to F$ une 
application bor\'elienne pour la topologie produit et continue le long
des horizontales $B\times \{t\}$ qui sera dite {\em
  de classe $BC^0$} (pour bor\'elienne et 
continue $C^0$). Pour tout $t\in T$, la restriction de $g$ d\'etermine
donc un \'el\'ement $g_\star(t)\in C(B,F)$, ce qui donne une
application $g_\star:T\to C(B,F)$ dite {\em verticale} de $g$. Le
r\'esultat suivant caract\'erise les applications de classe $BC^0$ 
comme celles dont les verticales sont bor\'eliennes. On remarquera que
la continuit\'e le long des horizontales est une condition
fondamentale. En effet, la verticale $g_\star$ d'une application
bor\'elienne quelconque $g:B\times T\to F$ est \`a valeurs dans
l'espace $\mathcal B(B,F)$ des applications bor\'eliennes de $B$ dans
$F$, qui n'est pas muni d'une structure mesurable naturelle.

\begin{prop}\label{prop:app-pri}
  Soient $B$ et $F$ deux espaces polonais, avec $B$ compact. Une
  application $g:B\times T\to F$ est de classe $BC^0$ si et 
  seulement si elle est continue le long des horizontales et
  l'application $g_\star:T\to C(B,F)$ est bor\'elienne.
\end{prop}
\begin{proof}
  Supposons que $g$ est de classe $BC^0$. On veut montrer que $g_\star$
  est bor\'elienne.  Soit $U_i$ une base d\'enombrable de $F$ et soit
  $U_i^*$ la base de $C(B,F)$ d\'efinie ci-dessus. Cette famille
  engendre alors la $\sigma$-alg\`ebre des bor\'eliens de $C(B,F)$; il
  suffit donc de montrer que l'ensemble 
  $$
  g_\star^{-1}(U_i^*)=\{t\in
  T~|~\forall x\in B, \,g(x,t)\in U_i\} 
  $$
  est un bor\'elien de $T$
  pour tout $i$. On remarque que 
  $$
  g_\star^{-1}(U_i^*)=\pi_T(B\times T-g^{-1}(F-U_i)) 
  $$
  o\`u $\pi_T$ est la projection de $B\times T$
  sur le deuxi\`eme facteur. Les fibres du bor\'elien $B\times T-
  g^{-1}(F-U_i)$ pour l'application $\pi_T$ sont ouvertes dans $B$,
  qui est par hypoth\`ese un espace polonais compact. En particulier
  elles sont r\'eunion d\'enombrable de compacts. L'ensemble
  $g_\star^{-1}(U_i^*)$ est alors bor\'elien d'apr\`es le th\'eor\`eme
  \ref{thm:thm_kallman}.

  Inversement supposons que $g_\star$ est bor\'elienne; on veut
  montrer qu'il en est de m\^eme pour $g$. On doit montrer que l'image
  inverse par $g$ de tout bor\'elien de $F$ est un bor\'elien de
  $B\times T$, et pour cela on peut se ramener aux ouverts d'une base
  d\'enombrable $U_i$ de $F$. Il s'agit donc de montrer que les
  ensembles 
  $$
  E_i=\{(x,t)\in K\times T~|~g(x,t)\in U_i\} 
  $$
  sont bor\'eliens dans $B\times T$.

  Soit $A$ un sous-ensemble de $B$. L'application de
  restriction 
  $$
  \mathfrak r_A:C(B,F)\to C(A,F) 
  $$
  est continue, donc bor\'elienne, relativement aux topologies
  compactes-ouvertes. En particulier l'ensemble 
  $$
  T_i(A)=\{t\in T~|~\forall x\in A,~g(x,t)\in U_i\}.  
  $$
  est bor\'elien pour tout $i$ et de plus $A\times T_i(A)\subset
  E_i$. Pour une base d\'enombrable $V_j$ de $B$ nous avons en
  particulier $\bigcup_j V_j\times T_i(V_j)\subset E_i$. 

  On compl\`ete
  la preuve de la proposition en d\'emontrant l'inclusion inverse
  $E_i\subset \bigcup_j V_j\times T_i(V_j)$. Pour ce faire nous prenons
  un point dans $E_i$ ou, ce qui est la m\^eme chose, une paire
  $(x,t)\in B\times T$ telle que $g(x,t)\in U_i$. Mais $g_\star(t)$
  \'etant continue, il existe un voisinage ouvert $V$ de $x$ telle que
  $g(\bar x,t)\subset U_i$ pour tout $\bar x\in V$. Puisque $V_j$ est
  une base de $B$, il existe un $j$ tel que $x\in V_j\subset V$. On a
  alors $g(\bar x,t)\in U_i$ pour tout $\bar x\in V_j$, ce qui
  signifie que $t\in T_i(V_j)$ et d\'emontre l'inclusion cherch\'ee.
\end{proof}

Le lecteur remarquera que l'on pourrait remplacer sans beaucoup de
difficult\'es l'hypoth\`ese de 
compacit\'e de $B$ par la compacit\'e locale. Mais nous n'aurons besoin
dans la suite que du cas compact.

\subsection{Approximation simpliciale}\label{cha:app-sim}
Soit $\H$ un complexe simplicial fini et $\L$ un complexe simplicial
connexe et localement fini. On note $|\H|$ et $|\L|$ les
r\'ealisations g\'eom\'etriques de $\H$ et $\L$ respectivement. 
Soit $g:|\H|\to |\L|$ une application continue. On rappelle qu'une
application simpliciale $h:|\H|\to |\L|$ est une {\em
  approximation simpliciale} de $g$ si 
$$
g(star(v,\H))\subset star(h(v),\L) 
$$
pour tout sommet $v\in \H^0$ o\`u $star(v,\H)$ d\'esigne l'\'etoile
ouverte de $v$ dans $\H$. On notera comme d'habitude $sd(\H)$ le
complexe simplicial obtenu par subdivision barycentrique de $\H$ et
$sd^n(\H)=sd(sd^n(\H))$ avec $sd^0(\H)=\H$. La preuve du r\'esultat
classique suivant peut \^etre trouv\'ee dans \cite{Mu1}:

\begin{thm}[d'approximation simpliciale finie] \label{thm:app-fin}
  Soit $\H$ un complexe simplicial fini et $f:|\H|\to |\L|$ une
  application continue. Alors il existe un $n\in \N$ et une
  approximation simpliciale $h:|sd^n(\H)|\to |\L|$ de $f$.
\end{thm}

On notera $\Sigma(|\H|,|\L|)$ l'ensemble des applications continues
$f\in C(|\H|,|\L|)$ qui sont simpliciales modulo une subdivision
barycentrique de $\H$. On munit $C(|\H|,|\L|)$ de la m\'etrique
$d^*_\L$ de la convergence uniforme relative \`a la m\'etrique
simpliciale $d_\L$. On appellera {\em application d'approximation
  simpliciale} toute application bor\'elienne 
$$
\mathfrak a:C(|\H|,|\L|)\to \Sigma(|\H|,|\L|) 
$$
telle que $\mathfrak a(g)$ est une
approximation simpliciale de $g$ pour toute $g\in C(|\H|,|\L|)$. Le
th\'eor\`eme \ref{thm:app-fin} montre qu'une application de ce type
existe, mais ne dit rien \`a propos de la r\'egularit\'e qu'on peut en
esp\'erer. Puisque $\H$ est un complexe fini,
l'ensemble $\Sigma(|H|,|\L|)$ est d\'enombrable, et on ne peut donc pas
esp\'erer qu'elle soit continue. Le r\'esultat
suivant montre qu'elle peut \^etre suppos\'ee bor\'elienne:

\medskip
\begin{thm}\label{thm:app-sim-2}
  Pour tout complexe simplicial fini $\H$ et tout complexe simplicial
  d\'enombrable $\L$ il existe une application d'approximation
  simpliciale bor\'elienne 
  $$
  \mathfrak a:C(|\H|,|\L|)\to \Sigma(|\H|,|\L|).  
  $$
\end{thm}
\begin{proof}
  On remarque d'abord que pour toute $g\in C(|\H|,|\L|)$ il existe un
  entier $n$ tel le diam\`etre de l'image par $g$ de l'\'etoile de
  tout sommet de $sd^n(\H)$ est inf\'erieur o\`u \'egal \`a
  $1/2$. Ceci est une cons\'equence d'une part de la 
  continuit\'e de $g$ et d'autre part de la compacit\'e de $|\H|$. On
  d\'efinit $\mathfrak n(g)$ comme \'etant le plus petit des ces
  entiers $n$. Il est tr\`es facile de voir que l'application 
  $$
  \mathfrak n:C(|\H|,|\L|)\to \N 
  $$
  est bor\'elienne. On posera $Q_m=\mathfrak
  n^{-1}(m)$ pour tout $m\in \N$.

  Pour toute $h\in \Sigma(|\H|,|\L|)$ soit $d(h)$ le plus
  petit entier $n$ tel que $h$ est simpliciale relativement \`a
  $sd^n(\H)$ et $\L$. On note $B(h)$ l'intersection du bor\'elien
  $Q_{d(h)}$ avec la boule de centre $h$ et rayon $1/2$. La famille
  des $B(h)$, o\`u $h$ parcourt les \'el\'ements de
  $\Sigma(|\H|,|\L|)$, est d'apr\`es \ref{thm:app-fin} un recouvrement
  d\'enombrable de $C(|\H|,|\L|)$. On num\'erote les \'el\'ements
  $\Sigma(|\H|,|\L|)$ et 
  on pose $$
  X_{k+1}=B(h_{k+1})-X_k $$
  avec $X_1=B(h_1)$.

  Pour toute $g\in C(|\H|,|\L|)$ on d\'efinit $\mathfrak
  a(g)$ comme \'etant \'egale \`a $h_k$ pour $x\in X_k$. L'application
  ainsi construite est bor\'elienne car constante sur les \'el\'ements
  d'une partition bor\'elienne d\'enombrable. Il ne reste qu'\`a
  montrer que $\mathfrak a(g)$ est une approximation simpliciale de
  $g$. Par construction $g$ est dans la boule de centre $\mathfrak
  a(g)$ et rayon $1/2$.  Ceci implique que pour tout sommet $v$ de
  $sd^{\mathfrak n(g)}(\H)$, la distance entre $\mathfrak a(g)(v)$ et
  $g(v)$ est $< 1/2$. Une application directe de l'in\'egalit\'e
  triangulaire montre alors que $$
  g(star(v,sd^{\mathfrak n(g)}(\H)))
  $$
  est dans la boule de centre $\mathfrak a(g)$ et rayon $1$. Pour
  conclure il suffit de remarquer que la boule de centre $w$ et rayon
  $1$ relative \`a la m\'etrique simpliciale est contenue dans
  l'\'etoile de $w$ pour tout sommet $w\in \L^0$.
\end{proof}

\section{Th\'eorie de l'obstruction}\label{sec:theor-de-lobstr}
Nous \'elaborons ici une th\'eorie de l'obstruction analogue \`a celle
de \cite{Eil,St,Hu} adapt\'ee au cadre des feuilletages.
C'est une th\'eorie g\'en\'erale qui r\'esout compl\`etement le
probl\`eme de l'extension d'un morphisme de feuilletages
transversalement  mesurable, ainsi que celui de la construction de
sections transversalement mesurables d'un fibr\'e sur $M$. Nous
remarquons au passage qu'une ``th\'eorie de l'obstruction
feuillet\'ee continue'' n'est pas envisageable. Par exemple
consid\'erons un plongement du tore $\mathbb S^1\times \mathbb S^1$
dans $\R^3$ de sorte que le $0$ est dans la composante connexe
born\'ee du compl\'ementaire de l'image. Par projection on obtient une
application surjective dans la sph\`ere $\mathbb S^2$ qui ne peut pas
\^etre \'etendue au tore plein $\mathbb D^2\times \mathbb S^1$. Par
contre, si on suppose le tore feuillet\'e par les horizontales, on
remarque que la restriction de cette application \`a toute feuille
$\mathbb S^1\times \{*\}$ est extensible au disque correspondant
$\mathbb D^2\times \{*\}$. En d'autres mots nous avons un cocycle
feuillet\'e continu qui est nul sans que l'application soit
contin\^ument extensible. 

Ici $F$ d\'esigne un espace localement compact triangulable et
$n$-simple, i.e. connexe par arcs et tel que l'action de $\pi_1(F,x)$
sur $\pi_n(F,x)$ est triviale quel que soit $x\in F$. On peut
d\'efinir alors le groupe $\pi_n(F)$ sans faire r\'ef\'erence au point
base $x\in F$.  Nous renvoyons le lecteur \`a \cite{St} pour plus de
pr\'ecisions \`a propos de cette d\'efinition. Nous remarquerons
simplement les deux faits suivants:
\begin{itemize}
\item Si $F$ est $1$-simple, alors $\pi_1(F)$ est un groupe ab\'elien.
\item Une application continue $f:\mathbb S^n\to F$ d\'etermine un
  \'el\'ement dans $[f]\in \pi_n(F)$ qui ne d\'epend que de
  l'application $f$ et de l'orientation choisie sur $\mathbb S^n$.
\end{itemize}

\subsection{Les lemmes fondamentaux}\label{sec:les-lemm-fond}
Pour tout espace m\'etrique compact $X$ on consid\`ere l'espace
des applications continues $C(X,F)$ muni de la topologie
compacte-ouverte. On note $\mathbb D^{n+1}$ le disque unit\'e ferm\'e 
dans $\R^n$ et $\mathbb S^n$ la sph\`ere unit\'e. On identifie le
disque $\mathbb D^n$ \`a la calotte sud de $\mathbb S^n$ via un
hom\'eomorphisme qu'on fixe une fois pour toutes. Nous avons alors
la suite d'inclusions 
$$
\xymatrix{ \mathbb D^n\ar[r]^{i_-} & \mathbb
  S^n\ar[r]^{i_+} &\mathbb D^{n+1}, } 
$$
qui induit par restrictions une suite
d'applications continues, donc bor\'eliennes,
$$
\xymatrix{ C(\mathbb D^{n+1},F)\ar[r]^{i_+^*} & C(\mathbb
  S^n,F)\ar[r]^{i_-^*} & C(\mathbb D^n,F) }. 
$$
Il est bien connu que l'application ${i_-^*}$ est
surjective, i.e. toute application continue sur la calotte sud
s'\'etend \`a toute la sph\`ere. Ce n'est pas le cas
de ${i_+^*}$. Son image est le sous-espace ferm\'e $C_0(\mathbb
S^{n-1},F)$ form\'e par les applications homotopes \`a z\'ero. 

Rappelons que le $n$-i\`eme groupe d'homotopie $\pi_n(F)$ est le
groupe engendr\'e par le semi-groupe $[\mathbb S^n,F]$ des classes
d'homotopie d'applications continues de $\mathbb S^n$ dans $F$ avec
la composition usuelle. On consid\`ere
l'application continue $[\cdot]:C(\mathbb S^n,F)\to 
\pi_n(F)$ qui assigne \`a toute fonction $g$ sa classe
d'homotopie $[g]$. La suite d'applications bor\'eliennes
\begin{equation}\label{eq:obs}
  \xymatrix{
    C(\mathbb D^{n+1},F)\ar[r]^{i^*_+} & C(\mathbb
    S^n,F)\ar[r]^{[\cdot]} & \pi_n(F)   
  }
\end{equation}
est exacte dans le sens o\`u $\mathrm{im\,} i^*_+=\ker\, [\cdot]$,
l'espace $\ker\, [\cdot]$ \'etant par d\'efinition $C_0(\mathbb
S^{n-1},F)$.

\begin{lem}\label{lem:fond-1}
  L'application $i^*_+$ admet une section bor\'elienne $$
  \mathfrak
  s:C_0(\mathbb S^n,F)\to C(\mathbb D^{n+1},F).  $$
\end{lem}
\begin{proof}
  On suppose $F$ et $\mathbb D^{n+1}$ munis d'une triangulation. Celle
  de $\mathbb D^{n+1}$ induit une triangulation sur son
  bord $\mathbb S^n$. On notera $\Sigma(*,F)$ l'espace de fonctions
  qui sont simpliciales apr\`es une subdivision barycentrique de
  l'espace de d\'epart. On commence par remarquer que la restriction
  de $i^*_+$ \`a l'espace $\Sigma(\mathbb D^{n+1},F)$ est \`a fibres
  d\'enombrables et son image est le bor\'elien $\Sigma(\mathbb
  S^n,F)\cap C_0(\mathbb S^n,F)$. Elle admet par le th\'eor\`eme
  \ref{thm:thm_kallman} une section bor\'elienne: 
  $$
  \mathfrak b:\Sigma(\mathbb S^n,F)\cap C_0(\mathbb
  S^n,F) \to \Sigma(\mathbb D^{n+1},F) 
  $$
  
  Pour la construction de la section $\mathfrak s$ on proc\`ede en
  trois \'etapes: 
  
  \medskip\noindent {\bf (I)} On consid\`ere
  $\mathfrak{h_a}:C(S^n,F)\to C(\mathbb S^n\times [0,1])$ d\'efinie
  par: $$
  \mathfrak{h_a}(f)(x,t)= tf(x)+(1-t)\mathfrak a(f)(x) $$
  o\`u
  $\mathfrak a(f)$ est l'approximation simpliciale de $f$ donn\'ee par
  le th\'eor\`eme \ref{thm:app-sim-2}. L'appli\-cation $\mathfrak{h_a}$
  est bor\'elienne comme combinaison lin\'eaire d'applications
  bor\'eliennes. 

  \medskip\noindent {\bf (II)} L'application $\mathfrak
  m:C(\mathbb S^n,F)\to \Sigma(\mathbb D^{n+1},F)$ d\'efinie par: 
  $$
  \mathfrak m(f)=\mathfrak{b\circ a}(f)
  $$
  est bor\'elienne comme composition de deux applications
  bor\'eliennes. 

  \medskip\noindent {\bf (III)} L'application $\mathfrak s(f)$
  d\'efinie par recollement $\mathfrak m(f)$ et
  $\mathfrak{h_a}(f)$ le long des sph\`eres $\mathbb S^n$ et $\mathbb
  S^n\times 1$ est une section bor\'elienne de $i^+_*$ par
  construction.
\end{proof}

Soient $a,b:\mathbb D^n\to F$ deux applications continues qui
co\"\i ncident sur le bord $\mathbb S^{n-1}$. On peut voir $a$ et $b$
comme des applications d\'efinies sur la calotte nord et la calotte
sud de la sph\`ere $\mathbb S^n$ respectivement et on note
$a*b:\mathbb S^n\to F$ le recollement de $a$ et $b$. Il est tr\`es
facile de voir que l'ensemble des paires $(a,b)$ d'applications
d\'efinies sur le disque et co\"\i ncidant sur le bord est un ferm\'e
$K$ de $C(\mathbb D^n, F)\times C(\mathbb D^n,F)$ et que le
recollement est une op\'eration continue sur $K$.

En munissant la sph\`ere $\mathbb S^n$ de l'orientation qui induit
l'orientation canonique sur la calotte sud (donc l'orientation
oppos\'ee sur la calotte nord), puis passant aux classes d'homotopie,
on d\'etermine un \'el\'ement $[a*b]\in \pi_n(F)$. Le choix de
l'orientation fait que $[a*b]=-[b*a]$. Plus g\'en\'eralement,
supposons que les restrictions de $a$ et $b$ \`a la sph\`ere $\mathbb
S^{n-1}$ ne sont pas identiques, mais seulement homotopes. Toute
homotopie $h:\mathbb S^{n-1}\times [0,1]\to F$ entre ces deux
restrictions s'\'etend donc aux deux bases de $\mathbb D^n\times
[0,1]$ de mani\`ere \'evidente. En identifiant le bord de $\mathbb
D^n\times [0,1]$ \`a la sph\`ere $\mathbb S^n$ par un hom\'eomorphisme
qui envoie $\mathbb S^{n-1}\times [0,1]$ sur la couronne de rayon $1/2$
centr\'ee dans l'\'equateur, nous avons une application 
$$
a*_h b :\mathbb S^n\to F 
$$
et en orientant la sph\`ere $\mathbb S^n$ comme
pr\'ec\'edemment, on obtient un \'el\'ement $[a*_h b]\in \pi_n(F)$ qui
d\'epend de $a$ et de $b$, mais aussi de l'homotopie $h$. Nous avons
bien s\^ur $[a*b]=[a*_{id}~b]$.

Consid\'erons maintenant deux applications
$\alpha,\beta:\partial\Delta_{n+1}\to F$ qui sont homotopes en
restriction au $(n-1)$-squelette de $\Delta_{n+1}$. On se fixe une
homotopie $h$ entre ces restrictions, pour chaque $n$-face orient\'ee
$\tau$ de $\Delta_{n+1}$ on construit l'\'el\'ement $[\alpha_\tau*_h
\beta_\tau]\in \pi_n(F)$, o\`u $\alpha_\tau$ et $\beta_\tau$ d\'esignent les
restrictions de $\alpha$ et $\beta$ \`a la face $\tau$. On obtient
alors de fa\c con tr\`es simple le r\'esultat suivant:

\begin{lem}\label{lem:fond-2}
  Pour $\alpha$, $\beta$ et $h$ comme ci dessus nous avons
  l'identit\'e:
  \begin{equation}
    \sum_\tau [\alpha_\tau*_h\beta_\tau]=[\alpha]-[\beta]
  \end{equation}
  o\`u $\tau$ parcourt les $n$-faces de $\Delta_{n+1}$ munies de
  l'orientation induite par celui-ci.
\end{lem}

La paire $(\mathbb S^n,\mathbb D^n)$ ayant le type l'homotopie de la
sph\`ere point\'ee $(\mathbb S^n,*)$, toute application continue $f\in
C(\mathbb D^n,F)$ peut \^etre \'etendue \`a toute la sph\`ere $\mathbb
S^n$ de sorte que la classe d'homotopie de l'extension soit un
\'el\'ement de $\pi_1(F)$ fix\'e \`a l'avance. Le r\'esultat suivant
montre que cette extension peut \^etre choisie de fa\c con
bor\'elienne: 

\begin{lem}\label{lem:fond-3}
  Pour tout $\alpha\in \pi_1(F)$ il existe une application
  bor\'elienne $$
  \mathfrak s_\alpha:C(\mathbb D^n,F)\to C(\mathbb
  D^n,F) $$
  telle que, pour toute $g\in C(\mathbb D^n,F)$, on a:
  \begin{enumerate}
  \item $g$ et $\mathfrak s_\alpha(g)$ co\"\i ncident sur $\mathbb
    S^{n-1}$;
  \item $[g*\mathfrak s_\alpha(g)]=\alpha$.
  \end{enumerate}
\end{lem}
\begin{proof}
  On identifie la sph\`ere $\mathbb S^n$ au
  bord du disque \`a coins $\mathbb D^n\times [0,1]$ par un
  hom\'eomorphisme qui envoie la calotte nord sur la base
  sup\'erieure $\mathbb D^n\times \{1\}$. On note $\Sigma_+(\mathbb
  S^n,F)$ l'espace des applications continues de $\mathbb S^n$ dans $F$
  qui sont simpliciales en restriction \`a cette calotte, i.e. sur la
  base sup\'erieure $\mathbb D^n\times \{1\}$,
  et $\Sigma_-^\alpha(\mathbb S^n,F)\subset \Sigma_+(\mathbb S^n,F)$ le
  ferm\'e des applications dont la classe d'homotopie est $\alpha$.

  Comme pour la preuve du lemme pr\'ec\'edent on construit la section
  $\mathfrak s_\alpha$ en trois \'etapes:
  
  \medskip\noindent {\bf (I)}  La restriction de $i_-^*$ \`a
  $\Sigma_+^\alpha(\mathbb S^n,F)$ est \`a fibres d\'enombrables et
  son image est l'espace $\Sigma^\partial(\mathbb D^n,F)$ des
  applications continues qui sont simpliciales en restriction au 
  bord du disque. Elle poss\`ede
  donc une section bor\'elienne $\mathfrak n:\Sigma^\partial(\mathbb
  D^n,F)\to\Sigma_-^\alpha(\mathbb S^n,F)$.

  \medskip\noindent {\bf (II)}  On consid\`ere alors l'application
  continue $\partial: C(\mathbb D^n,F)\to C(\mathbb S^{n-1},F)$
  d\'efinie par restriction au bord, puis $\mathfrak
  {h_a}:C(\mathbb S^{n-1},F)\to C(\mathbb S^{n-1}\times [0,1])$
  d\'efinie comme dans la preuve du lemme pr\'ec\'edent. Enfin on
  d\'efinit $f\in C(\mathbb D^n,F)\mapsto \mathfrak R(f)\in
  \Sigma^\partial(\mathbb D^n,F)$ par recollement de $\mathfrak
  {h_a}\circ \partial(f)$ et $f$.

  \medskip\noindent {\bf (III)} L'application $\mathfrak
  s_\alpha$ d\'efinie par $\mathfrak s_\alpha(f)=i^*_+\circ \mathfrak
  n \circ \mathfrak R(f)$ v\'erifie par construction les conditions
  requises.
\end{proof}

\subsection{Triangulations de feuilletages}\label{sec:triang}
Soit $(M,\mathcal F)$ une vari\'et\'e compacte feuillet\'ee. Nous
introduisons ici une notion de triangulation de feuilletages, similaire
\`a celle de \cite{HL}. Notre d\'efinition est n\'eanmoins bien plus
souple, en particulier une triangulation au sens de \cite{HL}
d\'efinit une triangulation selon notre d\'efinition.

\paragraph{Piles.}
On consid\`ere un espace m\'etrique compact connexe $\Omega$ et un
espace bor\'elien standard $T$. Une {\em pile} de $(X,\mathcal F)$ est
donn\'ee par un bor\'elien $\Pi\subset X$ et un isomorphisme
bor\'elien: 
$$
\pi:\Omega\times T\to\Pi $$
tel que, pour tout $t\in T$,
la restriction $\pi_t:\Omega\to M$ de $\pi$ \`a $\Omega\times \{t\}$
est un plongement de $\Omega$ 
dans une feuille de $\mathcal F$. Les \'el\'ements $\Omega$, $T$ et
$\pi$ seront appel\'es respectivement la {\em base}, la {\em
  verticale} et le {\em  param\'etrage} de la pile. Les {\em plaques}
de $(\Pi,\pi)$ sont les ensembles $\Pi_t=\pi(\Omega\times \{t\})$.

\paragraph{Triangulations.}\label{sec:tri}
Une {\em triangulation bor\'elienne (de classe $C^r$)} de $(M,\mathcal
F)$ est une famille $\K=\{\K_L|L\in\mathcal F\}$ de triangulations de
classe $C^r$ des feuilles de $\mathcal F$ telle que  pour chaque entier
$0\leq p\leq \dim \mathcal F$ il existe une quantit\'e d\'enom\-brable
de piles (de classe $BC^r$):
$$
\pi^p_i:\Delta_p\times T_i^p\to \Sigma^p_i\quad,\quad i\in I_p 
$$
de base $\Delta_p$ le $p$-simplexe standard et v\'erifiant les
propri\'et\'es suivantes:
\begin{enumerate}
\item les plaques de $\pi_i^p$ sont des $p$-simplexes de $\K$.
  
\item pour chaque $p$-simplexe $\sigma$ existe un seul $i\in I_p$ et
  un $t\in T_i$ tel que $$
  \sigma=\pi_i^p(\cdot,t)(\Delta_p).  $$
  On
  note alors $\pi_\sigma$ l'hom\'eomorphisme
  $\pi_i^p(\cdot,t):\Delta_p\to \sigma$.
\end{enumerate}

Par abus de langage on notera aussi $\K$ l'ensemble des simplexes de
la triangulation. C'est un complexe simplicial non connexe ni
s\'eparable mais \`a composantes connexes s\'eparables. On note
$\K^{(p)}$ 
l'ensemble des simplexes de dimension $p$ de $\K$, $\K^p$ l'ensemble
des simplexes de dimension $\leq p$ et $\K^{[p]}$ l'ensemble des
$p$-simplexes orient\'es, i.e. des paires form\'ees par un
$p$-simplexe plus une orientation de celui-ci. Pour $\mathcal H\subset
\K$ on note $|\mathcal H|\subset M$ la r\'eunion des simplexes de
$\mathcal H$. On a bien s\^ur $|\K|=M$.

Avec ces conditions l'ensemble $\K$ peut \^etre identifi\'e \`a
l'espace bor\'elien standard $\cup_{p,i} T_i^p$, les ensembles
$\K^{p}$ et $\K^{(p)}$ \'etant des bor\'eliens de $\K$. L'ensemble
des $p$-simplexes orient\'es $\K^{[p]}$ est muni \'egalement d'une
structure bor\'elienne standard qui fibre sur $\K^p$ avec une fibre \`a
deux points. On peut voir aussi $\K$ comme \'etant le bor\'elien de
$M$ form\'e par les barycentres de ses simplexes.

\paragraph{Cohomologie simpliciale.} \label{sec:cohom} Soit
$(M,\mathcal F)$ un feuilletage muni d'une triangulation bor\'elienne
$\K$. Soit $\Gamma$ un groupe ab\'elien bor\'elien, i.e. un groupe
ab\'elien muni d'une structure bor\'elienne standard pr\'eserv\'ee par
la somme et l'inversion. Une {\em $p$-cocha\^\i ne} de $\K$ est une
application bor\'elienne $c:\K^{[p]}\to \Gamma$ telle que
$c(-\tau)=-c(\tau)$, o\`u $-\tau$ est le simplexe $\tau$ muni de
l'orientation oppos\'ee. Le {\em cobord} de $c$ est la
$(p+1)$-cocha\^\i ne de $\K$ d\'efinie par $$
dc(\sigma)=\sum_{\tau\subset \sigma} c(\tau)\quad,\quad \sigma\in
\K^{[p+1]} $$
o\`u $\tau$ parcourt l'ensemble des $p$-simplexes de $\K$ contenus
dans $\sigma$ et munis de l'orientation induite par celui-ci. On
v\'erifie de la fa\,con usuelle que $d^2c=0$ pour toute cocha\^\i ne
$c$.

On note $C^p(\K,\Gamma)$ le groupe ab\'elien des $p$-cocha\^\i nes de
$\K$ et $C^*(\K,\Gamma)$ la somme directe de ces groupes. Puisque la
somme dans $\Gamma$ est une op\'eration bor\'elienne, le cobord d'une
cocha\^\i ne est bor\'elien. On a donc un op\'erateur nilpotent
$d:C^*(\K,\Gamma)\to C^*(\K,\Gamma)$ dont la cohomologie est not\'ee
$H^*(\K,\Gamma)$ et appel\'ee {\em cohomologie simpliciale de $\K$ \`a
  valeurs dans $\Gamma$}. On appelle comme d'habitude {\em
  $p$-cocycles} et {\em $p$-cobords} les $p$-cocha\^\i nes de $\K$
appartenant respectivement au noyau et \`a l'image de $d$.

\subsection{Le th\'eor\`eme central}\label{sec:theoreme-central}

Soit $(M,\mathcal F)$ un feuilletage muni d'une triangulation
mesurable $\K$. On se fixe un entier $p$ compris entre $0$ et $\dim
\mathcal F$ et on consid\`ere une application $g:|\K^p|\to F$ de classe
$BC^0$.

\begin{defn}\label{def:extensible}
  Soit $r>p$. On dira que $g$ est {\em $r$-extensible} s'il existe une
  application $\bar g:|\K^r|\to \F$ de classe $BC^0$ qui co\"\i ncide
  avec $g$ en restriction \`a $|\K^p|$.
\end{defn}

On va \'etudier dans cette section le probl\`eme de l'extensibilit\'e
d'une telle application $g$. Pour cela il est utile de
d\'ecouper le $p$-squelette en piles de $p$-simplexes, de sorte que
nous pouvons voir $g$ indiff\'eremment comme une application
$g^\star:\K^{p+1}\to C(\mathbb S^p,F)$ ou comme une application
$g_\star:\K^p\to C(\mathbb D^p,F)$, toutes les deux bor\'eliennes. Ces
deux points de vue sont \'equivalents en vertu de la proposition
\ref{prop:app-pri}.

Pour chaque $(p+1)$-simplexe orient\'e $\sigma\in \K^{[p+1]}$ on
consid\`ere: $$
c(g)(\sigma)=[g^*(\sigma)] \in \pi_p(F).  $$
On d\'efinit ainsi $(p+1)$-cocha\^\i ne $c(g)\in C^{p+1}(\K;\pi_p(F))$
qui s'av\`ere \^etre un cocycle (cf. \cite{St}) que l'on appelle
commun\'ement {\em cocycle d'obstruction}.

\medskip Le th\'eor\`eme suivant constitue le coeur de la th\'eorie de
l'obstruction feuillet\'ee:

\begin{thm}\label{thm:obst-central}
  Les trois propri\'et\'es suivantes sont v\'erifi\'ees:
  \begin{enumerate}
  \item Le cocycle $c(g)=0$ si et seulement si $g$ est
    $(n+1)$-extensible.
  \item Soient $g_0,g_1:|\K^p|\to F$ deux applications de classe
    $BC^0$ qui sont $BC^0$-homotopes en restriction \`a $|\K^{p-1}|$.
    Pour une telle homotopie $h$ il existe une $p$-cocha\^\i ne
    $\omega(g_0,h,g_1)\in C^p(\K,\pi_p(F))$ telle que: $$
    d\omega(g_0,h,g_1)=c(g_0)-c(g_1) $$
    En particulier
    $[c(g_0)]=[c(g_1)]\in H^{p+1}(\K,\pi_p(F))$. On appelle
    $\omega(g_0,g_1)$ la {\em cocha\^\i ne diff\'erence}.
  \item Soit $g_0:|\K^p|\to F$ une application de classe $BC^0$. Pour
    toute $p$-cocha\^\i ne $\omega\in C^p(\K;\pi_p(F))$ il existe une
    application $g_1:|\K^p|\to F$ de classe $BC^0$ qui co\"\i ncide
    avec $g_0$ sur $|\K^{p-1}|$ et telle que: $$
    \omega=\omega(g_0,g_1).  $$
  \end{enumerate}
\end{thm}
\begin{proof}
  Les trois conditions sont des corollaires plus ou moins directs des
  lemmes \ref{lem:fond-1}, \ref{lem:fond-2} et \ref{lem:fond-3}. On
  reprend donc les notations introduites au \S\ref{sec:les-lemm-fond}.

  \medskip \underline{Preuve de 1}: Il est clair que si $g$ est
  $(p+1)$-extensible alors $g^\star$ est \`a valeurs dans $C_0(S^p,F)$
  et $c(g)=0$. R\'eciproquement si $g^\star$ est \`a valeurs dans
  $C_0(S^p,F)$, on d\'efinit une extension $\tilde g$ de $g$ en
  posant $\tilde g_\star = \mathfrak s\circ g^\star$.

  \medskip \underline{Preuve de 2}: Soit $h^\star:\K^p\to C(\mathbb
  S^{p-1}\times [0,1],F)$ l'application bor\'elienne associ\'ee \`a
  l'homotopie $h$. La cocha\^\i ne $\omega(g_0,h,g_1)=[g_0^\star *_{h^\star}
  g_1^\star]$ v\'erifie la condition requise
  d'apr\`es le lemme \ref{lem:fond-2}.

  \medskip \underline{Preuve de 3}: Il suffit de poser
  $g_1^\star(\sigma)=\mathfrak s_\alpha \circ g_0^\star(\sigma)$ si
  $\omega(\sigma)=\alpha$. L'application $g_1^\star$ est bien
  bor\'elienne puisque $\pi_1(F)$ est d\'enombrable.
\end{proof}

Consid\'erons un deuxi\`eme feuilletage $(N,\mathcal G)$ muni
d'une triangulation bor\'elienne $\L$, et soit $\phi:M\to N$ une
application de classe $BC^0$ simpliciale le long de chaque
feuille. Elle d\'etermine 
de la fa\c con usuelle un homomorphisme de cocha\^\i nes
$\phi^\sharp:C^*(\L,\Gamma)\to C^*(\K,\Gamma)$ pour n'importe quel
groupe ab\'elien bor\'elien $\Gamma$. La preuve du r\'esultat suivant
est compl\`etement standard (voir par exemple \cite{St,Hu}):

\begin{thm}\label{thm:natural-cocycle}
  Soit $g:|\L^p|\to F$ une application de classe $BC^0$. Alors $g\circ
  \phi:|\K^p|\to F$ est de classe $BC^0$ et on a: $$
  c(g\circ
  \phi)=\phi^\sharp(c(g)).  $$
\end{thm}

\section{Classes caract\'eristiques
  feuillet\'ees}\label{sec:class-caract-feuill} 
Les d\'eveloppements de la section pr\'ec\'edente peuvent \^etre
g\'en\'eralis\'es au cas des fibr\'es localement triviaux sur $M$
ayant une fibre simpliciale. Nous suivrons dans la mesure du possible
la d\'emarche d\'ecrite par Steenrod \cite{St} suivant les travaux de
Eilenberg \cite{Eil}. On se fixe ici un feuilletage $(M,\mathcal F)$
muni d'une triangulation $\K$, et on consid\`ere un fibr\'e
topologique localement trivial $\xi=(E,p,M,F)$ dont la fibre $F$ est
un espace connexe localement compact triangulable et simple.

Avant de continuer on fixe quelques points de vocabulaire:

\begin{enumerate}
\item Une {\em $p$-section} de $\xi$ est une section de
  $\xi$ d\'efinie sur $|\K^p|$. On supposera que toutes les sections
  sont de classe $BC^0$. On dira qu'une telle section est {\em
    $r$-extensible} ($r>p$) s'il existe une $r$-section de classe
  $BC^0$ de $\xi$ qui co\"\i ncide avec $g$ sur le $p$-squelette.

\item Deux $p$-sections $g_0$ et $g_1$ de $\xi$ de classe $BC^0$ sont
  {\em homotopes} s'il existe une application $h: |\K^p|\times
  [0,1]\to E$ de classe $BC^0$ telle que:
  \begin{enumerate}
  \item $h(\cdot,t)$ est une $p$-section (a fortiori de classe $BC^0$
    de $\xi$ pour tout $t\in[0,1]$.
  \item $h(\cdot,0)=g_0$ et $h(\cdot,1)=g_1$.
  \end{enumerate}
  
\end{enumerate}

\subsection{Cohomologie simpliciale \`a valeurs dans un fibr\'e de
  coefficients.} \label{sec:cohom-coeff} Un fibr\'e de coefficients
sur $M$ est un fibr\'e principal $\mathbf \Gamma$ dont la fibre est un
groupe discret d\'enombrable $\Gamma$. Pour tout $x\in M$ on notera
$\Gamma_x$ la fibre de $\mathbf \Gamma$ au dessus de $x$; si $x$ est
le barycentre d'un simplexe $\sigma$ de $\K$ on posera
$\Gamma_x=\Gamma_\sigma$. 

\begin{ex}
  Soit $\xi$ un fibr\'e dont la fibre est un complexe simplicial
  simple $F$. On note $\mathbf \Pi_p(\xi)$ ou tout simplement $\mathbf
  \Pi_p$ le fibr\'e de coefficients associ\'e \`a $\xi$ ayant pour
  fibre le groupe $\pi_p(F)$. Il est d\'efini en rempla\c cant un
  cocycle de $\xi$, qui est \`a valeurs dans $Aut(F)$, par le cocycle
  \`a valeurs dans $Aut(\pi_p(F))$ obtenu par passage aux classes
  d'homotopie.    
\end{ex}

Une {\em $p$-cocha\^\i ne de $\K$ \`a valeurs dans $\mathbf \Gamma$}
est une application bor\'elienne $$
c:\K^{[p]}\to \mathbf \Gamma $$
telle que $c(\sigma)\in \Gamma_\sigma$ et $c(-\sigma)=-c(\sigma)$.
L'espace des $p$-cocha\^\i nes est un groupe ab\'elien dont la
structure est donn\'ee par celle de $\Gamma$ et que l'on notera
$C^{p}(\K;\mathbf \Gamma)$.

La trivialit\'e locale permet une identification entre
les groupes $\Gamma_\sigma$ et $\Gamma_\tau$ pour toute face $\tau$ de
$\sigma$. On peut ainsi d\'efinir le {\em cobord} d'une
$p$-cocha\^\i ne $c$ de $\K$ par la formule
\begin{equation}\label{eq:cob}
  dc(\sigma)=\sum_{\tau\subset \sigma} c(\tau)\quad,\quad \sigma\in \K^{[p+1]}
\end{equation}
o\`u $\tau$ parcourt les $p$-faces de $\sigma$ munies de l'orientation
induite par celle-ci. Puisque $\Gamma$ est discret et $\sigma$ est
simplement connexe, l'identification entre $\Gamma_\sigma$ et
$\Gamma_\tau$ est ind\'ependante de la trivialisation choisie. Le
cobord est donc un morphisme de groupes bien d\'efini
$d:C^p(\K,\mathbf \Gamma)\to C^{p+1}(\K,\mathbf \Gamma)$. On v\'erifie
aussi de la fa\c con usuelle que $d^2=0$; la cohomologie du morphisme
$d$  est appel\'ee {\em cohomologie simpliciale de $\K$ \`a valeurs
  dans $\mathbf \Gamma$} et not\'ee $H^*(\K,\mathbf \Gamma)$.

\subsection{Obstruction \`a l'extensibilit\'e des
  sections}\label{sec:class-caract} 
On se fixe un fibr\'e localement trivial $\xi=(E,p,M,F)$ au dessus de
$M$ comme ci-dessus et on consid\`ere $g:|\K^p|\to E$
une $p$-section de classe $BC^0$ de $\xi$. En se fixant une
trivialisation 
de $\xi$ au dessus de chaque $(p+1)$-simplexe $\sigma\in \K^{(p+1)}$
on identifie la fibre $F_x$ avec $F_\sigma$ pour tout $x\in \sigma$
o\`u $F_\sigma$ est la fibre de $\xi$ au dessus du barycentre de
$\sigma$, on peut voir la restriction de $g$ \`a $\partial\sigma$
comme une application $g^\star(\sigma)\in C(\mathbb S^p,F_\sigma)$. En
prenant la classe d'homotopie de $g^\star(\sigma)$ on d\'efinit une
$(p+1)$-cocha\^\i ne de $\K$ \`a valeurs dans le $p$-fibr\'e de
coefficients $\mathbf \Pi_p$ de $\xi$: $$
c(g):\sigma\in \K^{[p+1]}\to
[g^\star(\sigma)]\in\Pi_p.  $$
On remarque alors que puisque $F$ est
simple (en particulier $p$-simple), la classe d'homotopie de
$g^\star(\sigma)$ est ind\'ependante de la trivialisation choisie.
Nous avons dans ce cadre plus g\'en\'eral un \'equivalent du
th\'eor\`eme \ref{thm:obst-central}, avec une preuve tout \`a fait
analogue \`a celle d\'evelopp\'ee au \S~\ref{sec:theoreme-central}:

\begin{thm}\label{thm:obst-fibres}
  Les quatre propri\'et\'es suivantes sont v\'erifi\'ees:
  \begin{enumerate}
  \item Le cocycle $c(g)=0$ si et seulement si la section $g$ est
    $(n+1)$-extensible.
  \item Soient $g_0$ et $g_1$ sont deux $p$-sections de classe $BC^0$ de
    $\xi$ qui sont homotopes en restriction \`a $|\K^{p-1}|$. Pour une
    telle homotopie $h$ il existe une $p$-cocha\^\i ne
    $\omega(g_0,h,g_1)\in C^p(\K;\mathbf \Pi_p)$ telle que: 
    $$
    d\omega(g_0,h,g_1)=c(g_0)-c(g_1) 
    $$
    En particulier $[c(g_0)]=[c(g_1)]\in H^{p+1}(\K;\mathbf \Pi_p)$.
  \item Soit $g_0$ une $p$-section de classe $BC^0$ de $\xi$. Pour toute
    $p$-cocha\^\i ne bor\'elienne $\omega\in C^p(\K;\mathbf \Pi_p)$
    existe une $p$-section de classe $BC^0$ $g_1$ de $\xi$ de classe
    $BC^0$ qui co\"\i ncide avec $g_0$ sur $|\K^{p-1}|$ et telle que: 
    $$
    \omega=\omega(g_0,g_1).  
    $$
  \end{enumerate}
\end{thm}

Soit $(N,\mathcal G)$ un feuilletage muni d'une triangulation
bor\'elienne $\L$ et soit $\xi'=(E',p',N,F)$ un deuxi\`eme fibr\'e
topologique localement trivial de fibre $F$. Un
{\em morphisme bor\'elien} entre $\xi$ et $\xi'$ est une application
$\Phi:E\to E'$ de classe $BC^0$ qui envoie hom\'eomorphiquement fibre
sur fibre. Elle induit donc une application $\phi:M\to N$ entre les
bases d\'efinie par $p'\circ \Phi=\phi\circ p$. Le morphisme $\Phi$
est dit {\em simplicial} si l'application $\phi$ est simpliciale. Nous
avons:

\begin{thm}\label{thm:natural-cocyle-fibres}
  Soit $\Phi:\xi\to \xi'$ un morphisme simplicial et soit
  $g':|\L^p|\to E'$ une $p$-section bor\'elienne de $\xi'$. Alors il
  existe une et une seule $p$-section bor\'elienne $g$ de $\xi$ telle
  que $\Phi\circ g=g'\circ \phi$ et on a: $$
  \phi^\sharp(c(g'))=c(g).
  $$
\end{thm}

\subsection{Classe caract\'eristique
  feuillet\'ee}\label{sec:classe-caract-feuill} 
Le th\'eor\`eme \ref{thm:obst-fibres} permet la construction
topologique des classes caract\'eristiques feuillet\'ees d'un fibr\'e,
comme dans \cite{St}. On notera d\'esormais $\bf p$ le plus petit
entier pour lequel le groupe $\pi_{\bf p}(F)$ est non trivial.
 
\begin{lem}\label{lem:sections}
  Soit $\xi$ un fibr\'e localement trivial de fibre simple $F$ au
  dessus de $M$. Les deux propri\'et\'es suivantes sont v\'erifi\'ees:
  \begin{enumerate}
  \item $\xi$ poss\`ede une $({\bf p} -1)$-section de classe $BC^0$.
  \item Deux $({\bf p}-1)$-sections de classe $BC^0$ de
    $\xi$ sont homotopes.
  \end{enumerate}
\end{lem}
\begin{proof}
  L'existence de $({\bf p}-1)$-sections de classe $BC^0$ d\'ecoule de
  l'existence de $0$-sections, puis d'une r\'ecurrence sur $p$ en
  utilisant le th\'eor\`eme \ref{thm:obst-fibres}. Mais l'existence de
  $0$-sections de classe $BC^0$ est pratiquement triviale car,
  $|\K^0|$ \'etant une transversale bor\'elienne, la continuit\'e le
  long des feuilles est automatique, et il est tr\`es facile
  construire une section bor\'elienne d'un fibr\'e localement trivial
  au dessus d'un espace \`a base d\'enombrable, en l'occurrence $M$: on
  prend une suite $U_i$ d'ouverts trivialisants de $\xi$ et on
  d\'efinit des sections bor\'eliennes sur $U_0$,
  $U_1-U_0$, $U_2-(U_1\cup U_0)$ et ainsi de suite, qui se recollent
  entre elles en une section bor\'elienne de $\xi$, qu'on restreint
  ensuite \`a $|\K^0|$. 
  
  Pour montrer que deux $({\bf p}-1)$-sections $g_0$ et $g_1$ de $\xi$
  sont homotopes, on consid\`ere la vari\'et\'e \`a bord $M\times
  [0,1]$ muni du feuilletage $\mathcal F\times [0,1]$ dont les
  feuilles sont le produit de celles de $\mathcal F$ par
  l'intervalle $[0,1]$. La triangulation $\K$ d\'efinit des
  triangulations des feuilletages $(M,\mathcal F)\times 0$ et
  $(M,\mathcal F)\times 1$ qui s'\'etendent de mani\`ere \'evidente en
  une triangulation $\K'$ de tout le feuilletage $(M,\mathcal F)\times
  [0,1]$. On remarque alors que l'homotopie $h$ cherch\'ee n'est autre
  chose qu'une $p$-section du fibr\'e $\xi\times [0,1]$, pull-back de
  $\xi$ par la projection $M\times [0,1]\to M$. Plus pr\'ecis\'ement,
  c'est la restriction d'une $p$-section de $\xi\times [0,1]$ au
  sous-complexe $|\K^{{\bf p}-1}|\times [0,1]$ de $|\K'|$ qui
  co\"\i ncide avec $g_0$ et $g_1$ sur les bases $|\K^{{\bf
      p}-1}|\times 0$ et $|\K^{{\bf p}-1}|\times 1$. Elle existe
  d'apr\`es la propri\'et\'e 1 d\'emontr\'ee ci-dessus.
\end{proof}

Soient $g_0$ et $g_1$ deux $(\mathbf p-1)$-sections de classe $BC^0$ de
$\xi$. On consid\`ere une homotopie $h$ entre $g_0$ et $g_1$ comme
celle donn\'ee par le lemme pr\'ec\'edent. Soient et $\hat g_0$ et
$\hat g_1$ deux $\mathbf p$-extensions de $g_0$ et $g_1$
respectivement. Il existe en vertu du th\'eor\`eme
\ref{thm:obst-fibres}(2) une $\mathbf p$-cocha\^\i ne diff\'erence
$\omega(g_0,h,g_1)\in C^\mathbf p(\K,\mathbf \Pi_{\mathbf p})$ telle
que $d\omega(g_0,h,g_1)=c(g_0)-c(g_1)$. En particulier la classe de
cohomologie: 
$$
\mathfrak c(\xi,\K)=[c(\hat g_0)]=[c(\hat g_1)] 
$$
ne d\'epend ni des $(\mathbf p -1)$-sections ni des $\bf p$-extensions
choisies. C'est une classe canoniquement associ\'ee au fibr\'e $\xi$
et \`a la triangulation $\K$. Cette classe est appel\'ee {\em la
  classe caract\'eristique simpliciale} de $\xi$.

Nous avons le corollaire suivant du th\'eor\`eme \ref{thm:obst-fibres}:

\begin{thm}\label{thm:classe-zero-section}
  Soit $\xi$ un fibr\'e localement trivial de fibre simple $F$ et
  $(M,\mathcal F)$ un feuilletage de dimension $n$. Si $\pi_p(F)=0$
  pour tout $p\leq n-2$, alors la classe $\mathfrak c(\xi,\mathcal
  \K)=0$ si et seulement si le fibr\'e $\xi$ poss\`ede une section
  de classe $BC^0$. 
\end{thm}

\section{Preuve des th\'eor\`emes A, B, C et D}
\label{sec:poin-hopf-preuve}
Nous appliquons les r\'esultats d\'emontr\'es dans les sections
pr\'ec\'edents \`a la preuve des th\'eor\`emes annonc\'es dans
l'introduction.

\subsection{Mesures transverses invariantes}
\label{sec:mesure_transv}

Soit $(M,\mathcal F)$ un feuilletage sur une
vari\'et\'e compacte. Une {\em transversale bor\'elienne} de
$(M,\mathcal F)$ est un bor\'elien de $M$ qui rencontre toute feuille
le long d'un ferm\'e discret de cette feuille. Deux transversales $T$
et $S$ sont {\em isomorphes} s'il existe une transformation bijective
bi-bor\'elienne $\gamma:T\to S$ telle que $\gamma(x)$ est dans la
m\^eme feuille que $x$ pour tout $x\in T$. 

Une {\em mesure transverse} est une application
$\sigma$-additive $\mu$ qui assigne \`a chaque transversale
bor\'elienne $T$ de $(M,\mathcal F)$ un nombre $\mu(T)\in [0,\infty]$.
Une telle mesure est {\em invariante} si $\mu(T)=\mu(S)$ pour
toute paire de transversales bor\'eliennes isomorphes $T$ et $S$. Elle
sera dite {\em finie} si $\mu(T)<\infty$ pour toute transversale
bor\'elienne compacte $T$. On appellera dans la suite {\em feuilletage
  mesur\'e} tout feuilletage muni d'une
mesure transverse invariante finie $\mu$.

Le lemme suivant clarifie la signification de l'invariance d'une
mesure transverse. Sa preuve, compl\`etement \'el\'ementaire, est
laiss\'ee  au lecteur.

\begin{lem} \label{lem:mesure-invariante} Soit $T$ et $S$ deux
  transversales mesurables de $(M,\mathcal F)$ et $\alpha:T\to S$ une
  transformation mesurable telle que $f(x)\in L_x$ pour tout $x\in T$.
  Soit $\mu$ une mesure transverse invariante sur $(M,\mathcal F)$.
  Alors pour toute fonction $f\in L^1(T,\mu)$ on a: 
  $$
  \int_T f(x) d\mu(x) = \int_S \left(\sum_{x\in \alpha^{-1}(y)}
    f(x)\right)d\mu(y).   
  $$
\end{lem}

Un ensemble $A\subset
M$ est dit {\em $\mathcal F$-satur\'e} ou tout simplement {\em
  satur\'e}  s'il est r\'eunion de feuilles de $\mathcal F$. Un satur\'e $A$ est
{\em $\mu$-n\'egligeable} si toute transversale bor\'elienne $T\subset
A$ est de mesure nulle. 

\begin{defn}
  Un objet d\'efinit sur $(M,\mathcal F,\mu)$ sera dit {\em de classe
    $MC^0$} (pour mesurable et continu) s'il existe un bor\'elien
  satur\'e $\mu$-n\'egligeable $A\subset M$ tel que l'objet est de
  classe $BC^0$ (i.e. globalement bor\'elien et continu le long des
  feuilles) en restriction \`a $M-A$.
\end{defn}

On termine ce paragraphe en rappelant la notion de mesure ergodique.

\begin{defn}
  Une mesure transverse invariante $\mu$ sur $(M,\mathcal F)$ est dite
  {\em ergodique} si pour tout bor\'elien satur\'e $A$, soit lui soit
  son compl\'ementaire est $\mu$-n\'egligeable. Un feuilletage muni
  d'une mesure ergodique est appel\'e un feuilletage mesur\'e ergodique.
\end{defn}

\subsection{Triangulations et champs $\mu$-finis}
\label{sec:mu-finis}
Soit $(M,\mathcal F,\mu)$ un feuilletage mesur\'e. Une triangulation
$\K$ de $(M,\mathcal F)$ est dite {\em $\mu$-finie} si
$\mu(\K^{(p)})<\infty$ pour tout $p$. On peut trouver une preuve du
r\'esultat suivant dans \cite{HL}:

\begin{prop}\label{prop:exist-triang}
  Tout feuilletage mesur\'e $(M,\mathcal F,\mu)$ de classe
  $C^r$ poss\`ede une triangulation $\mu$-finie de classe $C^r$. 
\end{prop}

On appellera {\em champ tangent} toute section ${\bf X}$ du fibr\'e
tangent $T\mathcal F$ continue le long des feuilles. On
consid\'erera deux types de champs tangents:
\begin{enumerate}
\item[(i)] Ceux {\em de classe $BC^0$} qui sont globalement bor\'eliens
  entant qu'applications de $M$ dans $T\mathcal F$;
\item[(ii)] Ceux {\em de classe $MC^0$} qui sont bor\'eliens quitte \`a
  enlever un satur\'e $\mu$-n\'egligeable de $M$.
\end{enumerate}

On notera $O_{\bf X}$ l'ensemble des z\'eros de ${\bf X}$. Un champ
tangent ${\bf X}$ est dit {\em transverse} ou {\em \`a z\'eros
  isol\'es} si la trace de $O_{\bf X}$ sur $\mu$-presque toute feuille
est un ferm\'e discret de la feuille. On remarquera que le caract\`ere
ferm\'e d\'ecoule automatiquement de la continuit\'e de ${\bf x}$ le
long de la feuille. Si $\bf X$ est de classe $BC^0$ alors $O_{\bf X}$
est une transversale bor\'elienne, tandis qu'en classe $MC^0$
l'ensemble $O_{\bf X}$ est seulement mesurable, i.e. bor\'elien modulo
un sous-ensemble de mesure nulle. L'indice
local $ind_{\bf X}(x)$ est en tout cas bien d\'efini pour
$\mu$-presque tout $x\in O_{\bf X}$ et d\'etermine une fonction
mesurable sur $O_{\bf X}$.

\begin{defn}
  Un champ tangent transverse ${\bf x}$ est dit {\em $\mu$-fini} si
  $ind_{\bf X}\in L^1(O_{\bf X},\mu)$. Dans ce cas l'int\'egrale
  $\int_{O_{\bf X}} ind_{\bf X}\,d\mu$ est un nombre r\'eel que nous
  appellerons {\em indice moyen} du champ par rapport \`a $\mu$.
\end{defn}

\subsection{Un peu d'homologie}
\label{sec:homologie}

On se fixe une triangulation $\mu$-finie $\K$ sur le feuilletage
mesur\'e $(M,\mathcal F,\mu)$, ainsi qu'un isomorphisme bor\'elien
$\K\simeq [0,1]$ dont l'utilit\'e est de munir l'espace des simplexes
d'un ordre bor\'elien. Cet ordre permet de d\'efinir des applications
bord bor\'eliennes $\partial_i:\K^{[p+1]}\to \K^{[p]}$ telles que
$\partial_0(\sigma)$,
$\partial_1(\sigma),\dots,\partial_{p+1}(\sigma)$ sont les $p$-faces
de $\sigma$ \'ecrites en ordre croissant et munies des orientations
induites. 

On veut d\'efinir le complexe $(C_*(\K,\R),\partial)$ des {\em
  cha\^ines bor\'eliennes r\'eeles} de $\K$. Pour cela prenons un
$(p-1)$-simplexe orient\'e $\tau\in \K^{[p-1]}$. Puisque les feuilles
sont localement compactes, l'\'etoile de $\tau$ est finie, i.e. il
existe un nombre fini de $p$-simplexes orient\'es $\sigma$ tels que
$\tau=\partial_i\sigma$. Alors pour toute cocha\^ine $z\in
C^{p-1}(\K,\R)$ la somme finie
$$
\tilde z(\tau)=\sum_{\tau=\partial_i\sigma} z(\sigma).  
$$
d\'efinit un \'el\'ement $\tilde z\in C^{p-1}(\K,\R)$. On d\'efinit le
complexe des cha\^ines en posant $C_p(\K,\R)=C^p(\K,\R)$ pour tout
$p\in \N$ et $\partial z=\tilde z\in C_{p-1}(\K,\R)$ pour tout $z\in
C_p(\K,\R)$. L'homologie de ce complexe est appel\'ee l'{\em homologie
  bor\'elienne r\'eelle} de $\K$ et not\'ee $H_*(\K,\R)$.

\subsection{L'indice de Kronecker}
\label{sec:indice-kronecker}
\begin{defn}
  Soient $c\in C^p(\K;\R)$ et $z\in C_p(\K;\R)$ une $p$-cocha\^\i ne
  et une $p$-cha\^\i ne. On d\'efinit l'{\em indice de Kronecker
    global} ou {\em moyen} de $z$ et $c$ par: 
  $$
  \langle c,z\rangle_\mu =\frac{1}{2}\int_{\K^{[p]}} c(\sigma)\cdot
  z(\sigma)~d\mu(\sigma) 
  $$
\end{defn}
On introduit le facteur $\frac{1}{2}$ pour
compenser le fait que l'on int\`egre deux fois la valeur
$c(\sigma)\cdot z(\sigma)=c(-\sigma)\cdot z(-\sigma)$. Remarquons
que l'int\'egrale ci-dessus n'est pas toujours d\'efinie ni
finie. Une paire cocha\^ine-cha\^ine $(c,z)$ sera dite {\em
$\mu$-finie} si $\int |c\cdot z|d\mu<\infty$. Dans ce cas $\langle
c,z\rangle_\mu$ est un nombre r\'eel bien d\'efini.

\begin{prop}\label{thm:indice-kronecker}
  Soit $c\in C^{p-1}(\K;\R)$ et $z\in C_p(\K;\R)$ des (co)cha\^\i nes
  de $\K$. Si les paires $(dc,z)$ et $(c,\partial z)$ sont
  $\mu$-finies alors: 
  $$
  \langle dc,z\rangle_\mu=\langle c,\partial z\rangle_\mu. 
  $$
\end{prop}
\begin{proof}
  Soit $\partial_i:\K^{[p]}\to \K^{[p-1]}$ les applications face.
  Rappelons qu'elles sont mesurables. Nous avons: $$
  \langle
  dc,z\rangle_\mu=\int_{\K^{[p]}} dc(\sigma)\cdot
  z(\sigma)\,d\mu(\sigma)= \sum_i \int_{\K^{[p]}} c(\partial_i\sigma)\cdot
  z(\sigma)\, d\mu(\sigma) $$
  $$
  \langle c,\partial z\rangle_\mu
  =\int_{\K^{[p-1]}} c(\tau)\cdot \partial z(\tau)\,d\mu(\tau) =
  \sum_i \int_{\K^{[p-1]}} \sum_{\partial_i\sigma=\tau}c(\tau)\cdot
  z(\sigma)\,d\mu(\tau) $$
  En appliquant le lemme
  \ref{lem:mesure-invariante} \`a la fonction $f:\K^{[p]}\to \R$
  d\'efinie par $f(\sigma)=c(\partial_i\sigma)\cdot z(\sigma)$ et \`a
  la transformation $\alpha=\partial_i$, on obtient l'\'egalit\'e
  entre les termes respectifs des sommes \`a droite. Ceci d\'emontre
  la proposition.
\end{proof}

\subsection{Orientation et cycle fondamental}
\label{sec:cycle_fondam}
On note $\mathcal O(\mathcal F)$ le rev\^etement des orientations de
$\mathcal F$; c'est un rev\^etement \`a deux feuillets de $M$
feuillet\'e par les rev\^etements des orientations des feuilles. Une
{\em orientation bor\'elienne} de $\mathcal F$ est une
section $\mathfrak o$ de $\mathcal O(\mathcal F)$ de classe
$BC^0$. Le feuilletage est bor\'eliennement orientable s'il poss\`ede
une orientation bor\'elienne. Il est clair qu'une orientation au sens
classique du fibr\'e vectoriel $T\mathcal F$ d\'efinit une
orientation bor\'elienne de $\mathcal F$, autrement dit un feuilletage
orientable est bor\'eliennement orientable. 

Une orientation bor\'elienne $\mathfrak o$ de $(M,\mathcal F)$
d\'efinit une orientation sur chaque feuille, et en particulier sur
les $n$-simplexes d'une triangulation 
$\K$. On note $\bf 1=\bf 1_{\mathfrak o}$ la $n$-cha\^ine bor\'elienne
qui vaut $1$ sur les $n$-simplexes de $\K$ munis de l'orientation
$\mathfrak o$. Il est bien connu qu'il s'agit d'un cycle que nous
appellerons le {\em cycle fondamental} de $(M,\mathcal F)$ (relatif \`a
l'orientation $\mathfrak o$).

\subsection{Preuve du th\'eor\`eme A}\label{sec:preuve-du-theoreme-A}

\subsubsection*{\'Etape I} 
On d\'efinit la caract\'eristique d'Euler d'une une triangulation
$\mu$-finie $\K$ de $(M,\mathcal F)$  par:
$$
\chi(\K,\mu)=\sum_i
(-1)^i\mu(\K^{(i)})=\int_\K (-1)^{\dim\sigma}d\mu(\sigma).  
$$

Dans cette premi\`ere \'etape nous prouvons le lemme suivant:

\begin{lem}\label{lem:champ-caracteristique}
  Si $\K$ est une triangulation de classe $C^2$ de $(M,\mathcal
  F,\mu)$ alors il existe un champ tangent $\bf Z$ de classe $BC^0$,
  $\mu$-fini et \`a z\'eros isol\'es sur $(M,\mathcal F,\mu)$ tel que  
  $$ 
  \chi(\K,\mu)=\int_{O_{\bf Z}} ind_{\bf  Z}(x)\,d\mu(x).  
  $$ 
\end{lem}
\begin{proof}
  Soient $\K'=sd(\K)$ et $\K''=sd^2(\K)$ les premi\`ere et deuxi\`eme
  subdivisions barycentriques de $\K$. Tout sommet $v\in\K''^0$ est
  dans l'int\'erieur d'un seul simplexe de $\K$. On note $\eta(v)$ le
  barycentre de ce simplexe. Nous avons ainsi une application
  mesurable $\eta:\K''^0\to \K'^0$ qui engendre par lin\'earit\'e une
  application simpliciale 
  $$
  \eta:|\K''|\to |\K'|.  
  $$
  Pour tout $x\in
  |\K''|$ le point $\eta(x)$ est dans le m\^eme simplexe que $x$ de
  sorte que le chemin lin\'eaire orient\'e $c_x$ d'origine $x$ et
  extr\'emit\'e $\eta(x)$ est bien d\'efini. Ce chemin est de
  classe $C^1$ car la triangulation est de classe $C^2$. On
  note ${\bf Z}(x)$ le vecteur tangent \`a $c_x$ en 
  $x$. Ceci d\'efinit un champ de vecteurs tangents aux feuilles dont
  les singularit\'es correspondent aux barycentres des simplexes de
  $\K$, et on peut v\'erifier facilement
  l'identit\'e suivante:
  \begin{equation}\label{eq:cha-car}
    ind_{\bf Z}(\sigma)=(-1)^{\dim \sigma}
  \end{equation}
  o\`u $ind_{\bf Z}(\sigma)$ d\'esigne l'indice de ${\bf Z}$ au
  barycentre de $\sigma$. En effet ${\bf Z}$ pointe vers le barycentre
  de $\sigma$ en tout point de l'int\'erieur de $\sigma$. Par
  cons\'equent $-\bf Z$ pointe vers le barycentre de $\sigma$ aux
  points int\'erieurs des simplexes de $\K''$ transverses \`a
  $\sigma$. Le champ ${\bf Z}$ a donc pour vari\'et\'e stable au point
  singulier $\hat \sigma$ le simplexe $\sigma$. Mais il est bien connu
  que l'indice d'une telle singularit\'e est \'egale \`a la parit\'e
  de la dimension de sa vari\'et\'e stable. Ceci montre l'identit\'e
  (\ref{eq:cha-car}). On compl\`ete la preuve du lemme en int\'egrant
  l'identit\'e (\ref{eq:cha-car}) par rapport \`a $\mu$.
\end{proof}

\subsubsection*{\'Etape II}
Soit maintenant ${\bf X}$ un champ tangent $\mu$-fini et \`a z\'eros
isol\'es. Quitte \`a d\'eformer l\'eg\`erement la triangulation $\K$ on
peut supposer que $O_{\bf X}$ ne rencontre pas le $(n-1)$-squelette
$|\K^{n-1}|$. Ce champ d\'etermine donc par restriction une
$(n-1)$-section de $T^1\mathcal F$ que nous notons $g_{\bf X}$. Pour
tout $n$-simplexe $\sigma\in \K^{[n]}$ muni de l'orientation induite par
celle du feuilletage on a par d\'efinition: 
$$
c(g)(\sigma)=\sum_{x\in O_{\bf X}\cap \sigma} ind_{\bf X}(x).  
$$
En
int\'egrant sur $\K^{[n]}$ l'identit\'e ci-dessus on obtient: 
$$
\langle c(g_{\bf X}),\mathbf 1\rangle_\mu = \frac{1}{2}\int_{\K^{[n]}}
c(g_{\bf X}) (\sigma)\cdot \mathbf 1(\sigma)\, d\mu(\sigma)= \int_{O_{\bf X}} ind_{\bf X}(x)\,
d\mu(x) 
$$
Consid\'erons maintenant le champ $\bf Z$ de l'\'enonc\'e
du lemme \ref{lem:champ-caracteristique}. Les cocycles d'obstruction
$c(g_{\bf X})$ et $c(g_{\bf Z})$ sont cohomologues d'apr\`es le
th\'eor\`eme \ref{thm:obst-fibres}. En appliquant la proposition
\ref{thm:indice-kronecker} on obtient: $$
\chi(\K,\mu)=\langle c(g_{\bf Z}),\mathbf
1\rangle_\mu= \langle c(g_{\bf X}),\mathbf 1\rangle_\mu = \int_{O_{\bf
    X}} ind_{\bf X}(x)\, d\mu(x).  $$

\subsubsection*{\'Etape III}
On conclut la preuve du th\'eor\`eme par le lemme suivant:

\begin{lem}\label{lem:euler=euler-triang}
  Soit $\K$ une triangulation $\mu$-finie de $(M,\mathcal F)$. Alors:
  $$
  \chi(M,\mathcal F,\mu)=\chi(\K,\mu).  $$
\end{lem}
\begin{proof}
  La caract\'eristique d'Euler de $(M,\mathcal F,\mu)$ est par
  d\'efinition l'ac\-couplement de la classe d'Euler $e(T\mathcal F)\in
  H^n(M,\Z)$ avec la classe de Ruelle-Sullivan $[C_\mu]\in H_n(M,\R)$.
  Rappelons que la classe d'Euler $e(T\mathcal F)$ est le pull-back de
  la classe d'Euler universelle $e_n\in H^n(BO(n))$ par une
  application classifiante $f:M\to BO(n)$ du fibr\'e $T\mathcal F$. Il
  est bien connu que la classe $e_n$ est la classe caract\'eristique
  au sens du \S\ref{sec:classe-caract-feuill} du fibr\'e universel en
  $(n-1)$-sph\`eres: $$
  \mathbb S^{n-1}\to E^1(n)\to BO(n).  $$
  Quitte \`a prendre une approximation simpliciale de $f$ on a par la
  naturalit\'e du cocycle d'obstruction (th\'eor\`eme
  \ref{thm:natural-cocyle-fibres}) l'identit\'e: 
  $$
  f^*(e_n)=\mathfrak c(T^1\mathcal F,\K).
  $$ 
  La classe feuillet\'ee $\mathfrak c(T^1\mathcal F,\K)$ est
  repr\'esent\'ee par le cocycle d'obstruction $c({\bf X})\in
  C^n(\K,\Z)$ associ\'e \`a un champ tangent $\mu$-fini \`a z\'eros
  isol\'es $\bf X$. Par cons\'equent: 
  $$
  \chi(M,\mathcal F,\mu)=\langle e(T\mathcal F),[C_\mu]\rangle=\langle
  c({\bf X}),\mathbf 1\rangle_\mu=\chi(\K,\mu) 
  $$
  ce qui compl\`ete la preuve du lemme.
\end{proof}

\subsection{Preuve du th\'eor\`eme B}\label{sec:preuve-du-theoreme-B}
L'implication (2)$\Rightarrow$(1) est la plus facile; c'est un
corollaire du th\'eor\`eme A. En effet soit $\bf X$ un champ tangent
$\mu$-fini \`a z\'eros isol\'es non d\'eg\'en\'er\'es, la non
d\'eg\'en\'erescence d'un z\'ero signifie que son indice local est
$\pm 1$; c'est le cas par exemple du champ tangent $\bf Z$ construit
dans le lemme~\ref{lem:champ-caracteristique}. Si
$\mu(O_{\bf X})<\epsilon$ alors $|\chi(M,\mathcal F,\mu)|\leq
\mu(O_{\bf X})<\epsilon$. Si pour tout $\epsilon>0$ on peut trouver $\bf
X$ v\'erifiant cela, alors $\chi(M,\mathcal F,\mu)=0$.

Nous prouvons dans ce qui suit l'implication (1)$\Rightarrow$(2). On
se fixe donc un feuilletage mesur\'e ergodique $(M,\mathcal F,\mu)$
\`a caract\'eristique d'Euler nulle, ainsi qu'une triangulation
$\mu$-finie $\K$ de $(M,\mathcal F)$.

\subsubsection*{Quelques notations et une proposition technique}
\label{sec:notations}

On supposera que $(M,\mathcal F,\mu)$ est munie d'une orientation et
que les $n$-simplexes de $\K$ sont munis de l'orientation induite. De
cette fa\c con on peut voir toute $n$-cocha\^ine comme une fonction
bor\'elienne $c:\K^{(n)}\to \R$ d\'efinie sur 
les $n$-simplexes non orient\'es. Pour tout bor\'elien $A\subset
\K^{(n)}$ on notera $\mathbf 1_A$ la cocha\^ine qui vaut $1$ sur $A$ et
$0$ partout ailleurs. Remarquons par exemple que le cocycle
fondamental $\mathbf 1=\mathbf 1_{\K^{(n)}}$. On notera $||c||=\int
|c|d\mu$ la norme $L^1(\mu)$ de la cocha\^ine $c$.  

On commence par d\'emontrer la proposition suivante: 

\begin{prop}\label{prop:prop-technique}
  Soit $c\in C^n(\K;\Z)$ une $n$-cocycle dont toutes les valeurs non
  nulles sont $\pm 1$. Si $\langle c,\mathbf 1\rangle_\mu=0$ alors il
  existe une suite de $(n-1)$-cocha\^\i nes $\eta_r$ telle que 
  $$
  \lim_{r\to\infty}||c-d\eta_r|| =0.
  $$
\end{prop}

\subsubsection*{Preuve de la proposition}
\label{sec:preuve-lem-tec}
Un {\em chemin de $n$-simplexes} est une suite
  $\sigma_1\dots\sigma_k$ de $n$-simplexes telle que $\sigma_i$ et
  $\sigma_{i+1}$ ont une face principale commune pour tout $i$. Les
  $(n-1)$-faces $\sigma_i\cap\sigma_{i+1}$ d'un chemin seront
  suppos\'es munies de l'orientation induite par $\sigma_{i+1}$, qui
  est l'oppos\'ee de celle induite par $\sigma_i$. On
  remarque enfin que l'on peut voir l'ensemble des chemins comme un
  bor\'elien de $\cup_k(\K^{(n)})^k$.  

  \medskip
  Soit $c\in C^n(\K;\Z)$ v\'erifiant les conditions de
  l'\'enonc\'e. Notons $T_+(c)$ et $T_-(c)$ les bor\'eliens de $\K^{(n)}$
  sur lesquels $c$ prend respectivement les valeurs $+1$ et
  $-1$. Remarquons que 
  puisque $\langle c,\mathbf 1\rangle_\mu=0$ alors on a
  $\mu(T_+(c))=\mu(T_-(c))$. On peut supposer que $\mu(T_+(c))>0$ car
  dans le cas contraire il n'y a rien \`a prouver. 
  
  \begin{lem}\label{lem:chemins-ergodique}
    Pour $\mu$-presque tout $\sigma\in T_-(c)$ il existe un chemin de
    $n$-simplexes $\sigma_1\dots \sigma_n$ tel que $\sigma=\sigma_1$
    et $\sigma_n\in T_+(c)$.
  \end{lem}
  \begin{proof}
    Soient $A_+$ et $A_-$ les satur\'es de $M$ form\'es par les feuilles
    de $\mathcal F$ qui rencontrent respectivement $T_+(c)$ et
    $T_-(c)$. Puisque ces deux bor\'eliens sont de mesure positive et
    que la mesure $\mu$ est ergodique, le satur\'e $A_+\cap A_-$ est
    de mesure totale, ce qui implique le lemme.
  \end{proof}
  
 On suppose que $\K$ est d\'efinie par une famille de prismes
  $\pi^p_i$ ($i\in I_p$) comme au \S\ref{sec:triang}, et on consid\`ere
  l'application bor\'elienne naturelle $i:\K^{(n)}\to I_n$ qui assigne
  \`a chaque $n$-simplexe de $\K$ l'indice du prisme auquel il
  appartient. On dira qu'un chemin $\sigma_1\dots\sigma_k\in
  C^k$ est de type $\alpha\in I_n^k$ si 
  $i(\sigma_l)=\alpha_l$ pour tout $l$. 

  \begin{lem}
    Deux chemins de m\^eme type co\"incident ou sont disjoints.
  \end{lem}
  \begin{proof}
    On raisonne par r\'ecurrence sur la longueur des chemins. Le
    r\'esultat est \'evident pour les chemins de longueur
    $1$. Supposons qu'il est vrai pour les chemins de longueur
    $k-1$.  Il a alors deux cas:
    \begin{itemize}
    \item[a)] $\sigma_1\cap\sigma'_1\neq\varnothing$, donc
      $\sigma_1=\sigma'_1$. Alors puisque deux faces principales d'un
      simplexe sont d'intersection non vide, on a
      $\sigma_2\cap\sigma'_2\neq\varnothing$, donc
      $\sigma_2=\sigma'_2$ et on conclut par l'hypoth\`ese de
      r\'ecurrence.
    \item[b)] Le m\^eme raisonnement montre que $\sigma_1\neq\sigma'_1$
      implique $\sigma_2\neq\sigma'_2$ et on conclut comme
      pr\'ec\'edemment que les deux chemins sont disjoints. 
    \end{itemize}
  \end{proof}

  Pour tout $i\in I_n$ on d\'esigne $T_\pm^i(c)=T_\pm(c)\cap \pi_i^n$,
  le bor\'elien  des $n$-simplexes de $T_\pm$ de type $i$. On notera
  $D(c)$ l'ensemble de tous les chemins reliant un 
  $n$-simplexe de $T_-(c)$ \`a un $n$-simplexe de $T_+(c)$.
  Tout chemin de $n$-simplexes $\sigma_1\dots \sigma_k\in D(c)$
  d\'etermine de 
  fa\c con \'evidente une $(n-1)$-cocha\^ine \`a support dans les
  faces du chemin et dont le cobord est
  $\mathbf 1_{\sigma_k} - \mathbf 1_{\sigma_1}$. En vertu du
  lemme pr\'ec\'edent on peut recoller les cocha\^ines correspondant
  aux chemins de m\^eme type $\alpha=i_1\dots i_k$ pour obtenir une
  cocha\^ine bor\'elienne $\omega_\alpha^c$ dont
  le cobord est la $n$-cocha\^ine $\mathbf 1_{T_+^{i_k}(c)}-\mathbf
  1_{T_-^{i_1}(c)}$.

  On num\'erote de fa\c con arbitraire les types de chemins et on
  d\'esigne par $\alpha(c)=i_1(c)\dots i_k(c)$ le premier type $\alpha$ tel que
  $||d\omega_\alpha^c||>0$, et on pose:  
  $$
  \hat c=c - d\omega_{\alpha(c)}^c.
  $$

  \begin{lem}
    Le cocycle $\hat c$ v\'erifie les hypoth\`eses de la proposition
    \ref{prop:prop-technique}. De plus on a:
    \begin{enumerate}
    \item[(i)] $T_\pm(\hat c)\subset T_\pm(c)$;
    \item[(ii)] $\alpha(\hat c)>\alpha(c)$.
    \end{enumerate}
  \end{lem}
  \begin{proof}
    Par hypoth\`ese $c=\mathbf 1_{T_+(c)}-\mathbf 1_{T_-(c)}$. On aura alors:
    \begin{equation*}\begin{split}
    \hat c & =(\mathbf 1_{T_+(c)}-\mathbf 1_{T_-(c)})-
    (\mathbf 1_{T_+^{i_k(c)}(c)}-\mathbf 1_{T_-^{i_1(c)}(c)})\\
    & = \mathbf 1_{T_+(c)-T_+^{i_k(c)}(c)}-\mathbf
    1_{T_-(c)-T_-^{i_1(c)}(c)} 
    \end{split}\end{equation*}
    ce qui montre que les seules valeurs non nulles de $\hat c$ sont
    $\pm 1$ et qui \'etablit la propri\'et\'e (i). De plus d'apr\`es
    la  proposition~\ref{thm:indice-kronecker} on a $\langle \hat
    c,\mathbf 1\rangle_\mu =\langle c,\mathbf 1\rangle_\mu=0$. 

    Par ailleurs on a:
    $$
    T_-(\hat c)= T_-(c)-T_-^{i_1(c)}(c)\quad\textrm{et}\quad T_+(\hat c)=
    T_+(c)-T_+^{i_k(c)}(c)
    $$
    donc il n'existe aucun chemin de $D(\hat c)$ de type
    $\alpha(c)$. La cocha\^ine $\omega_{\alpha(c)}^{\hat c}$ v\'erifie alors
    $||d\omega_{\alpha(c)}^{\hat c}|| =0$, ce qui implique
    $\alpha(\hat c)>\alpha(c)$ par d\'efinition.
\end{proof}

D\'efinissons deux suites  par r\'ecurrence:
\begin{enumerate}
\item $c_0=c$ et $\eta_0=\omega^c_{\alpha(c)}$;
\item $c_{r+1}=\hat c_r=c-d\eta_r$ et
  $\eta_{r+1}=\eta_r+\omega^{c_r}_{\alpha(c_r)}$.  
\end{enumerate}

La proposition~~\ref{prop:prop-technique} revient donc \`a prouver:
\begin{lem}
$\lim_{r\to\infty} ||c_r||=0.$
\end{lem}
\begin{proof}
  On consid\`ere un chemin quelconque $\sigma_1\dots
  \sigma_k\in D(c)$ de type $\alpha$. Puisque
  $\alpha(c_{r+1})>\alpha(c_r)$ il existe un $r$ tel que
  $\alpha(c_r)>\alpha$. En particulier $c_r(\sigma_1)=c_r(\sigma)=0$.
  Par cons\'equent l'ensemble $\cap_rT_-(c_r)$ est compos\'e
  seulement des $\sigma\in T_-(c)$ dont il ne part aucun chemin
  arrivant sur $T_+(c)$. Sa mesure est nulle d'apr\`es le
  lemme~\ref{lem:chemins-ergodique}. Par cons\'equent: 
  $$
  \lim_{r\to \infty} ||c_r||=\lim_{r\to
    \infty} 2\mu(T_-(c_r)) =0. 
  $$
\end{proof}

\subsubsection*{Fin de la preuve du th\'eor\`eme B}
\label{theoreme_B}

Soit $\bf Z$ un champ tangent $\mu$-fini \'a z\'eros non
d\'eg\'en\'er\'es. Quitte \`a faire une subdivision de $\K$ puis
\`a d\'eformer un peu le champ $\bf Z$ on peut supposer que $T_0$ et
$T_1$ ne rencontrent pas le $(n-1)$-squelette de $\K$ et que tout
$n$-simplexe contient tout au plus un z\'ero de $\bf Z$. Le champ $\bf
Z$ d\'etermine alors une $(n-1)$-section $g_{\bf Z}$ du fibr\'e en
sph\`eres $T^1\mathcal F$ dont le cocycle d'obstruction $c(g_{\bf Z})$
ne prend que des valeurs $\pm 1$ et $0$.  De plus d'apr\`es le
th\'eor\`eme A on a: 
$$
0=\chi(M,\mathcal F,\mu)=\int_{O_{\bf Z}} ind_{\bf
  Z}(x)d\mu(x)=\langle c(g_{\bf Z}),\mathbf 1\rangle_\mu.
$$
Le cocycle $c=c(g_{\bf Z})$ v\'erifie alors les conditions du lemme
pr\'ec\'edent, et il est donc limite (pour la norme $L^1(\mu)$) d'une
suite de cobords $d\omega_r$. Mais d'apr\`es la propri\'et\'e (2) du
th\'eor\`eme~\ref{thm:obst-fibres} il existe pour chaque $\omega_r$
une $(n-1)$-section $g_r$ de
$T^1\mathcal F$ telle que $\omega_r$ est la cocha\^\i ne diff\'erence
$\omega(g_r, g_{\bf Z}))$. Les cocyles d'obstruction $c(g_r)$ sont \`a
valeurs dans $\{+1,0,-1\}$ et nous avons $\lim_{r\to\infty} ||c(g_r)||
=0$. En vertu du th\'eor\`eme~\ref{thm:obst-fibres} on peut \'etendre
$g_r$ en un champ tangent sans Z\'ero au dessus des simplexes o\`u
$c(g_r)$ s'annule. On peut l'\'etendre aussi par lin\'earit\'e au
dessus des autres $n$-simplexes, mais admettant une singularit\'e non
d\'eg\'en\'er\'ee au barycentre du simplexe. Le r\'esultat est une
suite de champs tangents $\mu$-finis $\mathbf X_r$ \`a z\'eros non
d\'eg\'en\'er\'es qui v\'erifie:
$$
\lim_{r\to\infty} \mu(O_{\mathbf X_r})\leq \int_{O_{\mathbf X_r}}
|ind_{\mathbf X_r}|d\mu=\lim_{r\to\infty} ||c(g_r)|| = 0
$$
ce qui compl\`ete la preuve du th\'eor\`eme B.

\subsection{Preuve du th\'eor\`eme C}
\label{sec:theoreme_C}
Un feuilletage mesur\'e $(M,\mathcal F,\mu)$ (muni d'une triangulation
$\mu$) est {\em hypercompact} s'il poss\`ede une {\em filtration
  compacte simpliciale}, i.e. une suite croissante de
bor\'eliens simpliciaux \`a feuilles compactes $B_n\subset M$ telle que
$M-\cup_n B_n$ est un satur\'e $\mu$-n\'egligeable. Le r\'esultat
suivant est un corollaire facile du th\'eor\`eme de
Connes-Feldman-Weiss. Dans \cite{BH} le lecteur pourra en trouver une
preuve:  

\begin{prop}
  Un feuilletage mesur\'e est hypercompact si et seulement si il est
  moyennable. 
\end{prop}

On se fixe un feuilletage mesur\'e ergodique orientable $(M,\mathcal
F,\mu)$  qui est moyennable et \`a caract\'eristique d'Euler nulle. On
suppose qu'il est muni d'une triangulation $\mu$-finie
$\K$, d'une filtration compacte  simpliciale $B_n$ comme ci-dessus, et
d'une $(n-1)$-section $\mathbf s$ de $T^1\mathcal F$ dont le cocycle
d'obstruction $c(\mathbf s)$ est \`a valeurs dans $\{+1,
0,-1\}$. D'apr\`es la preuve du th\'eor\`eme A nous avons $\langle
c(\mathbf s),\mathbf 1\rangle_\mu =0$.

\subsubsection*{Une nouvelle proposition technique}
\label{sec:lem_technique_2}

Par des techniques analogues \`a celles de la
proposition~\ref{prop:prop-technique} on a:

\begin{prop}\label{lem:lemme-technique_2}
  Soit $c\in C^n(\K;\Z)$ un $n$-cocycle dont toutes les valeurs non
  nulles sont $\pm 1$ et $B_n$ une filtration compacte simpliciale de
  $(M,\mathcal F?\mu)$. Si $\langle c,\bf 1\rangle_\mu=0$ alors il
  existe une suite de $(n-1)$-cocha\^\i nes $\eta_r$ telle
  que: 
  \begin{enumerate}
  \item[(i)] Pour tout $r\in \N$ les cocha\^\i nes $\eta_r$ et
    $\eta_{r+1}$ co\"\i ncident sur $B_r$; 
  \item[(ii)] $\lim_{r\to\infty}||c-d\eta_r|| =0$;
  \end{enumerate}
\end{prop}
\noindent{\em Esquisse de d\'emonstration.}
 On d\'efinit des cocha\^ines $\omega^c_{\alpha(c)}$ comme dans la preuve
 de la proposition~\ref{prop:prop-technique} mais en consid\'erant
 seulement les chemins $D_r(c)\subset D(c)$ qui sont contenus dans
 un $B_r$. On peut construire ainsi par une r\'ecurrence analogue une
 suite de cocha\^ines $\omega_{\alpha(c_{r+1})}^{c_{r+1}}$ \`a
 support dans $B_{r+1}-\textrm{support } \omega_{\alpha(c_r)}^{c_r}$,
 de sorte que la cocha\^ine 
 $$
 \eta_{r+1}=\sum_{i=1}^{r+1} \omega_{\alpha(c_r)}^{c_r}
 $$ 
 co\"incide avec $\eta_r$ sur $B_{r+1}$. Enfin la propri\'et\'e (ii)
 est prouv\'ee par le m\^eme argument en tenant compte que $M-\cup_r
 B_r$ est $\mu$-n\'egligeable.\qed

\subsubsection*{Fin de la preuve}
On se fixe ${\mathbf Z}$ et $g_{\mathbf Z}$ comme
dans la preuve du th\'eor\`eme B, et on applique le
lemme~\ref{lem:lemme-technique_2} au cocycle d'obstruction
$c(g_{\mathbf Z})$. Nous avons une suite de cocha\^ines $\eta_r$ qui
par la propri\'et\'e (2) convergent sur tous les simplexes de
$\cup_r B_r$, i.e. sur $\mu$-presque tout simplexe. Puisque les
cocycles $d\eta_r$ sont \`a valeurs $+1$, $0$ ou $-1$, leur norme
$L^1(\mu)$ est born\'ee par $\mu(\K^{(n)})$ qui est finie par
hypoth\`ese. Le th\'eor\`eme de la convergence domin\'ee de Lebesgue
implique alors:
$$
||d\eta-c(g_{\mathbf Z})||=\lim_{r\to \infty}||d\eta_r-c(g_{\mathbf Z})||=0 
$$
o\`u $\eta$ est la cocha\^\i ne limite des $\eta_r$. On a alors
$d\eta=c(g_{\mathbf Z})$ sur $\mu$-presque tout $n$-simplexe. Le
th\'eor\`eme~\ref{thm:classe-zero-section} garantit alors l'existence d'une section de
$T^1\mathcal F$ de classe $MC^0$.

\subsection{Preuve du th\'eor\`eme D}\label{sec:preuve-du-theoreme-C}
Comme dit dans l'introduction, il nous reste \`a prouver (1)$\Rightarrow$(3).

\subsubsection{Les nombres de Betti feuillet\'es}
\label{sec:nombres-Betti}

Soit $(M,\mathcal F,\mu)$ un feuilletage mesur\'e ergodique muni d'une
triangulation $\mu$-finie $\K$. Pour chaque $p$-simplexe $\sigma\in
\K^{(p)}$ on consid\`ere le groupe des $p$-cocha\^\i nes de carr\'e
int\'egrable $$
C^n_{(2)}(K_\sigma)=l^2(\K_\sigma^{[p]}) $$
de la
feuille $\K_\sigma$ de $\K$ contenant $\sigma$. La famille
$\{C^n_{(2)}(\K_\sigma)\}_{\sigma\in \K^0}$ est un champ d'espaces de
Hilbert au dessus de $\K^0$ au sens de \cite{Dix}.  L'int\'egrale de
Hilbert: 
$$
C^n_{(2)}(\K)=\int_{\K^0}^{\oplus} C^n_{(2)}(\K_\sigma)d\mu(\sigma) 
$$
s'identifie naturellement au
sous-espace des $p$-cocha\^\i nes mesurables \`a carr\'e int\'egrable
$C^p_{(2)}(\K,\mu)=L^2(\K^{[p]},\mu)$ (voir \cite{Dix,Con1,Gab2}). A
tout sous-espace ferm\'e $H\subset C^n_{(2)}(\K,\mu)$ on associe un
nombre positif
$$
\dim_\mu(H)\in [0,+\infty] 
$$
appel\'e la {\em dimension de Murray-von
  Neumann} de $H$ qui d\'epend de la mesure $\mu$ et qui v\'erifie les
propri\'et\'es naturelles suivantes:
\begin{enumerate}
\item[i)] $\dim_\mu(H)=0$ si et seulement si $H=0$;
\item[ii)] $\dim_\mu(H_1\oplus H_2)= \dim_\mu(H_1)+ \dim_\mu(H_2)$;
\item[iii)] $\dim_\mu(C^n_{(2)}(\K,\mu))=\mu(\K^{(n)})=c_n(\K,\mu)$.
\end{enumerate}

\begin{defn}
  On d\'efinit le {\em $n$-i\`eme nombre de Betti} de $\K$ relatif \`a
  $\mu$ par $$
  b_n(\K,\mu)=\dim_\mu(\mathcal H^n_{(2)}(\K)) $$
  o\`u
  $\mathcal H^n_{(2)}(\K)$ est l'espace des $n$-cocha\^\i nes
  harmoniques de $C^n_{(2)}(\K)$, i.e. le noyau de l'op\'erateur
  $d+d^*$.
\end{defn}

La preuve du r\'esultat suivant est \'el\'ementaire compte tenu des
propri\'et\'es (ii) et (iii) de $\dim_\mu$ et de la d\'ecomposition de
Hodge de l'espace $C^n_{(2)}(\K)$ par rapport \`a l'op\'erateur $d$:

\begin{prop}
  Si $\K$ est une triangulation $\mu$-finie alors on a: $$
  \chi(\K,\mu)=\sum_{n=0}^\infty (-1)^n b_n(\K,\mu).  $$
\end{prop}

\subsubsection{L'isomorphisme de de Rham-Hodge}
\label{sec:de_Rham-Hodge}

Soit $(L,g)$ une vari\'et\'e de Riemann \`a g\'eom\'etrie born\'ee et
soit $\K_L$ une triangulation {\em $g$-born\'ee} dans le sens o\`u  le
diam\`etre et le volume des simplexes maximaux de $\K_L$ sont born\'es
sup\'erieurement et inf\'erieurement. On suppose aussi que les
coordonn\'ees barycentriques sont born\'ees. Le r\'esultat suivant est
d\^u \`a Dodziuk: 
\begin{thm}[\cite{Dod}]
  Sous les conditions ci-dessus on a des isomorphismes d'espaces de
  Hilbert: $$
  \H^p(L,g)=\H^p_{(2)}(\K_L)\quad (p=0,\dots,n) $$
  o\`u
  $\H^p(L,g)$ est l'espace des $p$-formes harmoniques de carr\'e
  int\'egrable de $(L,g)$.
\end{thm}

Supposons maintenant que $(M,\mathcal F,\mu)$ est un feuilletage
mesur\'e ergodique de dimension deux. Une m\'etrique de Riemann sur
$(M,\mathcal F)$ est une famille $\{g_L|L\in\mathcal F\}$ de m\'etriques
de Riemann sur les feuilles de $\mathcal F$. Les m\'etriques seront
suppos\'ees, comme tous les objets jusqu'\`a maintenant, de classe
$BC^0$ ou $MC^0$. Une m\'etrique de Riemann est dite \`a
{\em g\'eom\'etrie born\'ee} si toutes les m\'etriques $g_L$ sont \`a
g\'eom\'etrie born\'ee. Les bornes peuvent d\'ependre de la feuille.

La preuve du r\'esultat suivant est tr\`es facile en dimension deux,
cas dont nous avons besoin. Le lecteur pourra en trouver une
esquisse dans \cite{BH} Il est tout de m\^eme valable en
dimension quelconque et la preuve peut \^etre faite en adaptant celle
de \cite{HL}.  
\begin{prop}
  Soit $(M,\mathcal F,\mu)$ un feuilletage mesur\'e muni d'une
  m\'etrique de Riemann $g$ de classe $MC^0$,  g\'eom\'etrie born\'ee
  et volume $\mu$-fini. Alors il existe une triangulation $\mu$-finie
  $\K$ de  $(M,\mathcal F)$ qui est $g$-born\'ee en restriction \`a
  chaque feuille. 
\end{prop}

\subsubsection{Fin de la preuve}
\label{sec:fin-preuve-D}

On se fixe une m\'etrique de Riemann $g$ de classe $MC^0$, \`a
g\'eom\'etrie born\'ee et de volume $\mu$-fini sur un feuilletage
mesur\'e de dimension deux $(M,\mathcal F,\mu)$, ainsi qu'une
triangulation $\K$ comme celle de la proposition pr\'ec\'edente.
Supposons enfin que $\chi(M,\mathcal F,\mu)=\chi(\K,\mu)=0$. Nous
devons conclure que $g$ est une m\'etrique de Riemann parabolique en
restriction \`a $\mu$-presque toute feuille de $\mathcal F$.

Puisque $\mu$ est ergodique, on a deux alternatives:
\begin{enumerate}
\item $\mu$ est concentr\'e sur une feuille compacte;
\item L'ensemble des feuilles compactes est $\mu$-n\'egligeable.
\end{enumerate}

Le premier cas est trivial car r\'eduit le feuilletage \`a la feuille
compacte en question, pour qui le th\'eor\`eme D
est bien connu. On suppose donc qu'on est dans le deuxi\`eme cas. Les
feuilles \'etant $\mu$-presque toutes non compactes, et compte tenu du
fait qu'une surface de Riemann non compacte \`a g\'eom\'etrie born\'ee
ne poss\`ede pas de $i$-formes harmoniques de carr\'e int\'egrable
pour $i=0,2$, nous avons: 
$$
0=\chi(M,\mathcal F,\mu)=-b_1(\K,\mu)=-\dim_\mu \H^1_{(2)}(\K) 
$$
Or l'int\'egrale de Hilbert: $$
\H^1_{(2)}(\K)=\int_{\K^0}^\oplus \H^1_{(2)}(\K_x)d\mu(x)
$$
est nulle si et seulement si $\mu$-presque tous les
$\H^1_{(2)}(\K_x)$ sont nuls. D'apr\`es le th\'eor\`eme de Dodziuk on
a $\H^1(L,g)=0$ pour $\mu$-presque toute feuille $L$. Mais la feuille
$L$ \'etant non compacte, son rev\^etement universel est
conform\'ement \'equivalent au plan euclidien $\C$ ou au plan
hyperbolique $\mathbb H$. Mais il est bien connu (voir par exemple
\cite{FK}) que dans le deuxi\`eme cas nous aurons $\H^1(L,g)\neq 0$.
Nous sommes donc dans le premier cas, ce qui compl\`ete la preuve du
th\'eor\`eme D.

\end{document}